\documentclass[11pt]{article}
\usepackage{amsmath,amssymb,  amsthm, stmaryrd, MnSymbol, hyperref}
\usepackage{amsfonts}
\usepackage{graphicx}
\usepackage{color}
\usepackage{enumerate}

\newtheorem{theo}{Theorem}[section]
\newtheorem{fact}[theo]{Fact}
\newtheorem{coro}[theo]{Corollary}
\newtheorem{pro}[theo]{Proposition}
\newtheorem{lemma}[theo]{Lemma}
\newtheorem{re}[theo]{Remark}

\theoremstyle{definition}
\newtheorem{definition}[theo]{Definition}

\begin{document}
\title{\bf Dynamic Probability Logics: \\ Axiomatization \& Definability}
\author{Somayeh Chopoghloo {\footnote{\noindent School of Mathematics, Institute for Research in Fundamental Sciences (IPM), PO Box 19395-5746, Tehran, Iran.
			E-mail: s.chopoghloo@ipm.ir.}}
		\hspace{8mm}
		\vspace*{5mm}
		Massoud Pourmahdian {\footnote{\noindent  Corresponding author,
				School of Mathematics, Institute for Research in Fundamental Sciences (IPM),
				PO Box 19395-5746, Tehran, Iran; Department of Mathematics and Computer Science, Amirkabir University of Technology (Tehran Polytechnic), 15194, Tehran, Iran.
				E-mail: pourmahd@ipm.ir.}}	   	
}
\date{}
\maketitle
\begin{abstract}
Our aim in this paper is twofold. We first study probabilistic dynamical systems from logical perspective. To this purpose, we introduce the finitary {\em dynamic probability logic} ($\mathsf{DPL}$), as well as its infinitary extension $\mathsf{DPL}_{\omega_1}\!$. Both these logics extend the {\em  (modal) probability logic} ($\mathsf{PL}$) by adding a temporal-like operator $\bigcirc$ (denoted as dynamic operator) which describes the dynamic part of the system. We subsequently provide Hilbert-style axiomatizations for both $\mathsf{DPL}$ and $\mathsf{DPL}_{\omega_1}\!$. We show that while the proposed axiomatization for $\mathsf{DPL}$ is strongly complete, the axiomatization for the infinitary counterpart supplies strong completeness for each countable fragment $\mathbb{A}$ of $\mathsf{DPL}_{\omega_1}\!$.

Secondly, our research focuses on the (frame) definability of important properties of probabilistic dynamical systems such as measure-preserving, ergodicity and mixing within $\mathsf{DPL}$ and $\mathsf{DPL}_{\omega_1}$. Furthermore, we consider the infinitary probability logic $\mathsf{InPL}_{\omega_1}$ (probability logic with initial probability distribution) by disregarding the dynamic operator. This logic studies {\em Markov processes with initial distribution}, i.e. mathematical structures of the form $\langle \Omega, \mathcal{A}, T, \pi\rangle$ where $\langle \Omega, \mathcal{A}\rangle$ is a measurable space, $T: \Omega\times \mathcal{A}\to [0, 1]$ is a Markov kernel and $\pi: \mathcal{A}\to [0, 1]$ is a $\sigma$-additive probability measure. We verify that the expressive power of $\mathsf{InPL}_{\omega_1}$ is strong enough for showing that the {\em $n$-step transition probability $T^n$} of Markov kernel $T$, is $\mathsf{InPL}_{\omega_1}$-definable. This fact implies that many natural stochastic properties of Markov processes such as stationary, invariance, irreducibility and recurrence are $\mathsf{InPL}_{\omega_1}$-definable. These facts justify the importance of studying these particular extensions of $\mathsf{PL}$.
\end{abstract}
\textbf{MSC 2020:} 03B48, 03B70, 03C75, 60J05, 37H
\\
\\
\textbf{Keywords:} Modal probability logic, infinitary logic, axiomatization, completeness, definability, Markov processes, random dynamical systems.
\section{Introduction}
Our main concern in this article is to investigate (discrete-time) Markov stochastic processes augmented by a dynamic mapping from the modal logic point of view. These mathematical structures, which are called {\em dynamic Markov processes}, are of the form $\langle{\Omega, \mathcal{A}, T, f}\rangle$ where $\langle{\Omega, \mathcal{A}}\rangle$ is a measurable space, $T: \Omega\times \mathcal{A}\to [0, 1]$ is a Markov kernel and $f:\Omega\to \Omega$ is a measurable function.

In a somewhat broader context, the notion of {\em random dynamical systems} \cite{Arn} covers one of the most important classes of dynamical systems with probabilistic features. Typically, these systems contain stochastic processes, e.g. Markov processes possibly augmented by some additional dynamic structures that describe the dynamic behavior of the system. In a sense, our investigations lay in logical descriptions of certain special cases of random dynamical systems. These structures have diverse applications, from stochastic differential equations to finance and economics \cite{Bha}.

There are various logical approaches towards modeling of probability structures, among which we highlight first-order logic and propositional modal logic. All of them have advantages and possibly some restrictions. In the first-order approach, probability structures are ordinary first-order structures supplemented by probability measures. The first-order probability logic contains infinitely many quantifiers of the form $Px\geq r$, for each $r \in \mathbb{Q}\cap [0, 1]$. For a (probability) formula $\varphi(x)$, the formula $(Px\geq r)\varphi(x)$ means that `the probability of the set $\{x\;|\; \varphi(x)\}$ is at least $r$'. Keisler in \cite{Keisler2} denotes the finitary first-order probability logic by $\mathcal{L}_{\omega P}$, while its infinitary extension (with countably infinite conjunctions and disjunctions) by $\mathcal{L}_{\omega_1 P}$. This line of research is mostly devoted to providing axiomatizations and proving relevant soundness and completeness. In \cite{Keisler2}, there is a comprehensive survey of these results in which one could find the relevant research articles related to this subject. Further, Fajardo and Keisler in \cite{Fajardo} develop certain model theory for stochastic processes.

On the other hand, propositional modal logic provides a realm for reasoning about probability with a somewhat restricted expressive power compared to first-order logic, although decidable. In this approach, bounds on probability are treated as modal operators. So there are countably many probability modal operators $L_r$, for each $r\in \mathbb{Q} \cap [0,1]$. As the first-order case, the formula $L_r \varphi$ is interpreted by `the probability of $\varphi$ is at least $r$'. The resulting modal probability logic is denoted by  $\mathsf{PL}$.
Despite its weaker expressive power, there are advantages for which one can hold to  $\mathsf{PL}$.  Firstly, it is shown that this logic is decidable \cite{Zhou2009}, and therefore, in terms of computability, this logic is nice.  Secondly, the semantic machinery of modal logic provides a natural setting for considering Markov processes as semantical structures of this logic. This feature makes $\mathsf{PL}$ more favorable than $\mathcal{L}_{\omega P}$. Lastly, modal logic, as opposed to first-order logic,  supplies the notion of (frame) validity that can be employed to express global properties of probability structures.

There are numerous papers in this area dealing with axiomatization which demonstrate a number of completeness for  $\mathsf{PL}$ \cite{Aum1999, fagin, Gold2010, Hei, Kozen2, Zhou2009} and prove some nice semantical properties \cite{Deshar, Pour}. There is also infinitary version of $\mathsf{PL}$ denoted by $\mathsf{PL}_{\omega_1}$ \cite{Bara, Iko, Meier}.

Here we aim to introduce {\em dynamic probability logic} ($\mathsf{DPL}$), which is obtained by adding a temporal-like modal operator $\bigcirc$ (denoted as dynamic operator) to the language of probability logic (Section \ref{DPL}). Calling this as a dynamic operator stems in particular modal logic called {\em dynamic topological logic} \cite{Kremer1} which studies structures of the form $\langle{X, \tau, f} \rangle$ where $\langle{X, \tau} \rangle$ is a topological space and $f: X\to X$ is a continuous function. So the mathematical structures in this logic include, but are not limited to, abstract dynamical system in the sense of ergodic theory. It is worth nothing that there are a good deal of papers which study the combination of temporal and probability logics among which we notice \cite{Dod, Ogn, Shak}. In these contexts, the modal temporal operator $\bigcirc$ when it applies to a formula $\varphi$ is interpreted as `$\varphi$ holds in the next instance of time'. This interpretation, although related, is different from what we expect from the operator $\bigcirc$ whose interpretation expresses the iteration of a dynamic transformation.

Our contributions in this paper are divided into two parts. The first direction of our research is devoted to proposing a Hilbert-style axiomatization for $\mathsf{DPL}$ (Subsection \ref{The axiomatization of DPL}) and demonstrating its strong completeness for the class of all dynamic Markov processes based on standard Borel spaces (Theorem \ref{scompletef}). To this end, we use a canonical model construction based on special maximal finitely consistent subsets of formulas  called {\em saturated sets}. This approach is inspired by the proof of strong completeness for Markovian logics in \cite{Kozen2}. We further examine the logics of some important subclasses of dynamic Markov processes, including the class of all dynamic Markov processes of the form $\langle{\Omega, \mathcal{A}, T, f}\rangle$ that are {\em measure-preserving}, i.e., $T(w, f^{-1}(A))=T(f(w), A)$ for each $w\in \Omega$ and $a\in \mathcal{A}$  (Theorem \ref{MP and Pure}). We also present a logic for the class of all {\em abstract dynamical systems}, i.e. structures of the form $\langle{\Omega, \mathcal{A}, \mu, f}\rangle$ where $\langle{\Omega, \mathcal{A}, \mu}\rangle$ is a probability space and $f: \Omega\to \Omega$ is a measure-preserving function (Theorem \ref{PS}).

Our ideas naturally extend to introducing the {\em infinitary dynamic probability logic} (Section \ref{IDPL}). This logic, which is denoted by $\mathsf{DPL}_{\omega_1}$, allows countable conjunctions and disjunctions.
The expressive power of $\mathsf{DPL}_{\omega_1}$ is compatible with $\sigma$-additivity of probability measures. So within this logic,  many properties of probability can be naturally axiomatized, and hence, it is not hard to extend ideas from \cite{Bara, Iko} to show that there exists a weakly complete Hilbert-style axiomatization for this logic (Subsection \ref{The axiomatization of IDPL}).  Meanwhile, we show that whenever the logic is restricted to its countable fragments, the proposed axiomatization is strongly complete for the class of all dynamic Markov processes  (Theorem \ref{com2}). We should point out that while the canonical model introduced for the proof of strong completeness for each countable fragment $\mathbb{A}$ of $\mathsf{PL}_{\omega_1}$ in \cite[Subsection 5.2]{Iko} depends on $\mathbb{A}$, here we show that the canonical model of $\mathsf{DPL}$  in Theorem \ref{scompletef} can be served uniformly as a canonical model for each countable fragment of $\mathsf{DPL}_{\omega_1}$.

Our second contribution in this article is allocated to investigating (frame) definability of natural properties of dynamic Markov processes (Section \ref{Definability}).  We show that some dynamic properties such as measure-preserving, ergodicity and mixing are definable within $\mathsf{DPL}$ and $\mathsf{DPL}_{\omega_1}$. Moreover, we consider the {\em infinitary probability logic with initial distribution} ($\mathsf{InPL}_{\omega_1}$) by disregarding the dynamic operator. This logic studies {\em Markov processes with initial distribution}, i.e.  structures of the form $\langle\Omega, \mathcal{A}, T, \pi\rangle$ where $\langle\Omega, \mathcal{A}, T\rangle$ is a Markov process and $\pi: \mathcal{A}\to [0, 1]$ is a $\sigma$-additive probability measure (Subsection \ref{nstep}). We show that the strong expressive power of ${\mathsf{InPL}_{\omega_1}}$ would allow us to define  {\em $n$-step transition probabilities $T^n$} of Markov kernel $T$ (Lemma \ref{semequ}). From this, we conclude that many natural stochastic properties of Markov processes such as stationary, invariance, irreducibility and recurrence can be stated within $\mathsf{InPL}_{\omega_1}$. These results particularly show that $\mathsf{DPL}$ as well as $\mathsf{DPL}_{\omega_1}$ are natural and important extensions of $\mathsf{PL}$.
\section{Dynamic probability logic} \label{DPL}
In this section, we introduce the  dynamic probability logic and provide a strongly complete Hilbert-style axiomatization for it. Moreover, we investigate the logic of some important subclasses of dynamic Markov processes.
\subsection{Syntax and semantics} \label{Syntax and semantics 1}
Let $\mathbb{P}$ be a countable set of propositional variables. The language $\mathcal{L}_\mathsf{DPL}$ of {\em dynamic probability logic} is recursively defined by the following grammar:
\begin{align*}
	\varphi ::&=\ \ p\mid \neg\varphi \mid \varphi \land \varphi \mid L_r \varphi \mid \bigcirc\varphi
\end{align*}
where $p \in \mathbb{P}$ and $r\in \mathbb{Q} \cap [0,1]$. We employ the Boolean connectives $\top$, $\perp$, $\to$, $\leftrightarrow$ and $\vee$ in a usual way. For a formula $\varphi $, $L_r \varphi$ is interpreted as `the probability of $\varphi$ is at least $r$'. Moreover, $M_r \varphi$ is considered as an abbreviation for $L_{1-r} \neg\varphi$ and interpreted as `the probability of $\varphi $ is at most $r$'. Likewise, the additional probability operator $L_{r_1\dots r_k} \varphi$ is introduced to abbreviate $L_{r_1} \dots L_{r_k} \varphi$ where $r_1, \dots, r_k\in \mathbb{Q} \cap [0,1]$.
On the other hand, the modality operator $\bigcirc$ is read as the {\em dynamic} (or {\em next}) {\em operator} and has dynamic (or temporal) interpretation; $\bigcirc\varphi$ states that `$\varphi$ holds at the next instance of time'.

Next, we define various notions of {\em dynamic Markov processes}. These notions resemble the idea of  frames in modal logic. When endowed with valuation functions for propositional variables, they are considered as {\em dynamic Markov models}. Before that, let us first recall the definition of a Markov process from, e.g., \cite{Douc}.

\begin{definition} \label{Markov}
Let $\langle{\Omega, \mathcal{A}}\rangle$ be a measurable space. A function  $T: \Omega\times \mathcal{A}\to [0, 1]$  is said to be a {\em Markov kernel} (or {\em transition probability}) if it satisfies the following conditions:
\begin{itemize}
\item for each $w \in \Omega$, $T(w, .): A \mapsto T(w, A)$ is a ($\sigma$-additive) probability measure on $\mathcal{A}$;
\item  for each $A \in \mathcal{A}$,  $T(., A): w \mapsto T(w, A)$  is a measurable function on $\Omega$.
\end{itemize}
In this case, the triple $\langle{\Omega, \mathcal{A}, T}\rangle$ is called a {\em Markov process}  (or {\em Markov chain}) on the state space $\langle{\Omega, \mathcal{A}}\rangle$. In the following, we also write $T(w,A)$ as $T(w)(A)$.
\end{definition}

\begin{definition}\label{frames}
	\begin{itemize}
		\item[1.] A quadruple $\mathfrak{P} = \langle{\Omega, \mathcal{A}, T, f}\rangle$ is called a {\em dynamic Markov process} 
	 whenever $\langle{\Omega, \mathcal{A}, T}\rangle$ is a Markov process and  $f:\Omega\to \Omega$ is a measurable function. An element $w\in\Omega$ is called a {\em world} (or {\em state}).
		\item[2.] A dynamic Markov process  is said to be {\em measure-preserving} if for each $w\in \Omega$ and $A\in \mathcal{A}$, we have $T(w, f^{-1}(A)) = T(f(w), A)$.
		\item[3.] We say that a dynamic Markov process  is a {\em purely probabilistic} if $f=id_{\Omega}$. In these processes, the dynamic part is essentially discarded.
		\item[4.] Further, a dynamic Markov process  is a {\em dynamic probability space} if for each $w, v\in \Omega$ and $A\in \mathcal{A}$, we have $T(w, A) = T(v, A)$. That means there exists a single probability measure $\mu: \mathcal{A} \to [0, 1]$ such that  $T(w, A)= \mu(A)$ for each $w\in \Omega$ and $A\in \mathcal{A}$.  In this situation, we  represent this dynamic probability space as $\langle{\Omega, \mathcal{A}, \mu, f}\rangle$.
		\item[5.] Finally,  a measure-preserving dynamic probability space  is said to be an {\em abstract dynamical system}.
	\end{itemize}
\end{definition}

\begin{definition}
	A {\em dynamic Markov model}, or simply a {\em  model}, is a tuple $ \mathfrak{M} = \langle{\Omega, \mathcal{A}, T, f, v}\rangle$ where $\langle{\Omega, \mathcal{A}, T, f}\rangle$ is a dynamic Markov process  and  $v: \mathbb{P} \to \mathcal{A}$ is a valuation function which assigns to every propositional variable $p \in \mathbb{P}$ a measurable set $v(p)\in\mathcal{A}$.
\end{definition}

Let $ \mathfrak{M} = \langle{\Omega, \mathcal{A}, T, f, v}\rangle$ be a model. If
the quadruple $\langle{\Omega, \mathcal{A}, T, f}\rangle$ belongs to either of the special classes introduced  in Definition \ref{frames} then we denote this model by its underlying process. For example, a model  $ \mathfrak{M} $  is {\em measure-preserving} if its underlying process is so.

The satisfiability relation for arbitrary $\mathcal{L}_\mathsf{DPL}$-formulas in a given model $ \mathfrak{M} = \langle{\Omega, \mathcal{A}, T, f, v}\rangle$, is
defined inductively as follows:
\begin{itemize}
	\item[] $\mathfrak{M}, w \vDash p \;$ iff $\; w \in v(p)$,
	\item[] $\mathfrak{M}, w \vDash \neg\varphi \;$ iff $\; \mathfrak{M}, w \nvDash \varphi$,
	\item[] $\mathfrak{M}, w \vDash \varphi\land \psi\;$ iff $\; \mathfrak{M}, w \vDash \varphi$ and $\mathfrak{M}, w \vDash \psi$,
	\item[] $\mathfrak{M}, w \vDash L_r \varphi \;$ iff $\; T(w, [\![\varphi]\!]_{\mathfrak{M}}) \geq r$. Here
	$[\![\varphi]\!]_{\mathfrak{M}} = \{ w \in \Omega\;| \; \mathfrak{M}, w \vDash \varphi \}$,
	\item[] $\mathfrak{M}, w \vDash \bigcirc\varphi \;$ iff $\; \mathfrak{M}, f(w) \vDash \varphi$.
\end{itemize}
A set $\Gamma$ of $\mathcal{L}_\mathsf{DPL}$-formulas holds in the world $w$ if $\mathfrak{M}, w \vDash \varphi$ for each $\varphi \in \Gamma$.

\begin{re} \label{remeas}
	For any formula $\varphi$ of $\mathcal{L}_\mathsf{DPL}$ and model $\mathfrak{M}$, the measurability of $[\![\varphi]\!]_{\mathfrak{M}}$ yields easily from the following facts.
	\begin{align*}
	[\![p]\!]_{\mathfrak{M}} & \in \mathcal{A},\\
		[\![\neg\varphi]\!]_{\mathfrak{M}} &= \Omega \setminus [\![\varphi]\!]_{\mathfrak{M}},\\
		[\![\varphi\land \psi]\!]_{\mathfrak{M}} &= [\![\varphi]\!]_{\mathfrak{M}} \cap [\![\psi]\!]_{\mathfrak{M}},\\
		[\![L_r \varphi]\!]_{\mathfrak{M}} &= \{w\in\Omega\;|\; T(w, [\![\varphi]\!]_{\mathfrak{M}})\geq r\},\\
		[\![\bigcirc\varphi]\!]_{\mathfrak{M}} &= f^{-1} ([\![\varphi]\!]_{\mathfrak{M}}).
	\end{align*}
\end{re}

The following fact is needed for the proof of strong soundness theorem (Theorem \ref{sound1}). An algebraic form of this fact  is appeared in \cite{Kozen2013}.

\begin{fact} \label{continuity}
	Let $\varphi$ be a formula of $\mathcal{L}_\mathsf{DPL}$, $\mathfrak{M}$ be a model and $r_1, \dots, r_k \in \mathbb{Q} \cap [0,1]$. Then the set $[\![L_{r_1\dots r_k} \varphi]\!]_{\mathfrak{M}}$ is measurable. Furthermore, for each $r\in\mathbb{Q}\cap [0,1]$ we have
	$$[\![L_{r_1\dots r_k r} \varphi]\!]_{\mathfrak{M}} =\bigcap_{s<r} [\![L_{r_1\dots r_k s} \varphi]\!]_{\mathfrak{M}}.$$
\end{fact}
\begin{proof}
	The first part of the claim is clear. The second part, can be proved by induction on $k$.
	For the basic step when $k=0$, we have
	\begin{align*}
		[\![L_r \varphi]\!]_{\mathfrak{M}} &= \{ w\in\Omega\;|\; T(w, [\![\varphi]\!]_{\mathfrak{M}}) \geq r \}\\
		&= \{ w\in\Omega\;|\; \forall s<r, \;T(w, [\![\varphi]\!]_{\mathfrak{M}}) \geq s\}\\
		&= \bigcap_{s<r} \{ w\in\Omega\;|\;T(w, [\![\varphi]\!]_{\mathfrak{M}}) \geq s\}\\
		&= \bigcap_{s<r} [\![L_s \varphi]\!]_{\mathfrak{M}}.
	\end{align*}
	Now, for the induction step, suppose that $t \in \mathbb{Q} \cap [0,1]$ and $[\![\psi^r]\!]_{\mathfrak{M}} = \bigcap_{s<r} [\![\psi^s]\!]_{\mathfrak{M}}$ where $\psi^r = L_{r_1\dots r_k r} \varphi$ and $\psi^s = L_{r_1\dots r_k s} \varphi$. Then,
	\begin{align*}
		[\![L_t \psi^r]\!]_{\mathfrak{M}} &= \{ w\in\Omega\;|\; T(w, [\![\psi^r]\!]_{\mathfrak{M}}) \geq t \}\\
		&= \{ w\in\Omega\;|\;T(w, \bigcap_{s<r} [\![\psi^s]\!]_{\mathfrak{M}}) \geq t\}\\
		&= \{ w\in\Omega\;|\; \inf_{s<r} \big(T(w, [\![\psi^s]\!]_{\mathfrak{M}})\big) \geq t\}\\
		&= \{ w\in\Omega\;|\; \forall s<r, \;T(w, [\![\psi^s]\!]_{\mathfrak{M}})\geq t\}\\
		&= \bigcap_{s<r}\{ w\in\Omega\;|\;T(w, [\![\psi^s]\!]_{\mathfrak{M}}) \geq t\}\\
		&= \bigcap_{s<r}[\![L_t \psi^s]\!]_{\mathfrak{M}}.
	\end{align*}
\end{proof}

\begin{definition}
	Customarily, a formula $\varphi$ is {\em valid in a model} $\mathfrak{M}$, denoted by $\mathfrak{M}\vDash \varphi$, if $\mathfrak{M}, w \vDash \varphi$ for all $w \in \Omega$. Also, $\varphi$ is \textit{valid in a dynamic Markov process} $\mathfrak{P}$, denoted by $\mathfrak{P}\vDash \varphi$, if it is valid in every model based on $\mathfrak{P}$.
	Likewise, $\varphi$ is \textit{valid in a class $\mathcal{C}$} of processes, denoted by ${\vDash}_{\mathcal{C}}\; \varphi$, if it is valid in every element of $\mathcal{C}$.
	Further, 
	$\varphi$ is  a {\em semantic consequence} of a set $\Gamma$ over a class $\mathcal{C}$  if for every model $\mathfrak{M}$ based on an element in $\mathcal{C}$, if $\Gamma$ holds in a world $w$ in $\mathfrak{M}$, then so is $\varphi$. In this situation, we write $\Gamma \;{\vDash}_{\mathcal{C}} \;\varphi$. When  $\mathcal{C}$ is the class $\mathcal{DMP}$ of all dynamic Markov processes, we may write ${\vDash}\; \varphi$ and $\Gamma \;{\vDash} \;\varphi$ instead of ${\vDash}_{\mathcal{DMP}}\; \varphi$ and $\Gamma \;{\vDash}_{\mathcal{DMP}} \;\varphi$, respectively.
\end{definition}

\begin{definition}
	By {\em dynamic probability logic}, or $\mathsf{DPL}$, we mean the set of all valid formulas of $\mathcal{L}_\mathsf{DPL}$ over the class of all dynamic Markov processes.
\end{definition}
\subsection{The axiomatization of $\mathsf{DPL}$} \label{The axiomatization of DPL}
In this subsection, we present a Hilbert-style axiomatization for $\mathsf{DPL}$, denoted by $\mathcal{H}_{\mathsf{DPL}}$, and prove its strong soundness with respect to the class of all dynamic Markov processes.

\begin{definition}
The proof system $\mathcal{H}_{\mathsf{DPL}}$ has the following axiom schemes:
\begin{align*}
& \textnormal{Taut}&&\text{All propositional tautologies }\\
&\textnormal{FA$_1$}&&L_0 \perp &&\text{(finite additivity)}\\
&\textnormal{FA$_2$}&&L_r \neg\varphi \to \neg L_s\varphi, \text{ if } r+s > 1 &&\text{}\\
&\textnormal{FA$_3$}&&L_r (\varphi\land \psi) \land L_s(\varphi\land \neg\psi) \to L_{r+s}\varphi, \text{ if }  r+s \leq 1&&\text{}\\
&\textnormal{FA$_4$}&& \neg L_r (\varphi\land \psi) \land \neg L_s(\varphi\land \neg\psi) \to \neg L_{r+s}\varphi, \text{ if } r+s \leq 1&&\text{}\\
&\textnormal{Mono}&&L_1 (\varphi\to\psi) \to (L_r\varphi\to L_r \psi)  &&\text{(monotonicity)}\\
&\textnormal{Func$_\bigcirc$}&&\bigcirc\neg\varphi\leftrightarrow \neg\bigcirc\varphi &&\text{(functionality)}\\
&\textnormal{Conj$_\bigcirc$}&&\bigcirc(\varphi \wedge \psi) \leftrightarrow (\bigcirc\varphi\wedge\bigcirc\psi) &&\text{(conjunction)}
\end{align*}
and the inference rules:
\begin{align*}
&\textnormal{MP}&&\frac{\displaystyle{\;\varphi\to\psi\;\;\;\varphi\;}}{\displaystyle{\psi}}&& \text{(modus ponens)}\;\\
& \textnormal{GArch$_{\bigcirc^n, r}$} && \frac{\displaystyle{\;\{\psi\to{\bigcirc}^n L_{r_1\dots r_k s}\varphi \;|\; s< r\}\;}}{\displaystyle{\psi \to {\bigcirc}^n L_{r_1\dots r_k r} \varphi}}, \text{ for }  n\in \mathbb{N} && \text{(generalized Archimedean)} \\
&\textnormal{Nec$_{L_1}$}&&\frac{\displaystyle{\;\varphi\;}}{\displaystyle{L_1 \varphi}} && \text{(necessitation)}\\
&\textnormal{Nec$_\bigcirc$}&&\frac{\displaystyle{{\;\varphi\;}}}{\displaystyle{\bigcirc\varphi}} && \text{}
\end{align*}
\end{definition}

\begin{re}
\begin{itemize}
\item[1.] The generalized Archimedean rule has {\em infinitary} nature in the sense that it allows us to derive a conclusion from a countably infinite set of premises.
\item[2.]  GArch$_{\bigcirc^n, r}$ is a generalization of the following Arch rule proposed by Zhou in \cite{Zhou2009} for modal probability logic ($\mathsf{PL}$). In the presence of this rule, the weak completeness theorem for $\mathsf{PL}$ is derived \cite[Theorem 3.19]{Zhou2009}.
\begin{align*}
&\textnormal{Arch} && \frac{\displaystyle{\;\{\psi\to L_s\varphi \;|\; s< r\}\;}}{\displaystyle{\psi \to L_r \varphi}}&& \text{\;\;\;\;\;\;\;\;\;\;\;\;\;\;\;\;\;\;\;\;\;\;\;\;\;}\;\;\;\;\;\;\;\;\;\;\;\;\;\;\;\;\;\;\;\;\;\;\;\;\;\;\;\;\;\;\;\;\;\;\;\;\;\;
\end{align*}
Furthermore, to have strong completeness for $\mathsf{PL}$, the following version of the generalized Archimedean rule was introduced by Kozen {\it et al.} in \cite{Kozen2}. This rule implies countable additivity in dynamic Markov processes.
\begin{align*}
&\textnormal{GArch} && \frac{\displaystyle{\;\{\psi\to L_{r_1\dots r_k s}\varphi \;|\; s< r\}\;}}{\displaystyle{\psi \to L_{r_1\dots r_k r} \varphi}}&& \text{\;\;\;\;\;\;\;\;\;\;\;\;\;\;\;\;\;\;\;\;\;\;\;\;\;}\;\;\;\;\;\;\;\;\;\;\;\;\;\;\;\;\;\;\;\;\;\;\;\;\;\;\;\;\;\;\;\;\;\;\;\;\;\;
\end{align*}
\end{itemize}
\end{re}

\begin{definition}\label{Theorem}
\begin{itemize}
\item[1.] We say that a formula $\varphi$ is a {\em theorem} of $\mathcal{H}_{\mathsf{DPL}}$, or is {\em derivable} in $\mathcal{H}_{\mathsf{DPL}}$, (denoted by $\vdash\varphi$) whenever there exists a sequence $\varphi_0, \varphi_1, \dots, \varphi_{\alpha+1}$ ($\alpha$ is a finite or infinite countable ordinal) of the formulas of $\mathcal{L}_\mathsf{DPL}$ with  $\varphi_{\alpha+1} =\varphi$ and for each $\beta \leq{\alpha+1}$,  $\varphi_\beta$ is either an axiom scheme or derived by applying one of the inference rules on some (possibly infinite) preceding formulas of the sequence.
\item[2.] A formula $\varphi$ is {\em derivable} from a set $\Gamma$ of assumptions in $\mathcal{H}_{\mathsf{DPL}}$ (denoted by $\Gamma  \vdash\varphi$), if there exists a sequence $\varphi_0, \varphi_1, \dots, \varphi_{\alpha+1}$ ($\alpha$ is a finite or infinite countable ordinal) of the formulas of $\mathcal{L}_\mathsf{DPL}$ with $\varphi_{\alpha+1} = \varphi$ and for each $\beta \leq{\alpha+1}$, either $\varphi_\beta$  is a member of $\Gamma$, or a theorem of $\mathcal{H}_{\mathsf{DPL}}$, or derived by applying one of the inference rules, other than Nec$_{L_1}$ and Nec$_\bigcirc$, on some preceding formulas of the sequence.
\item[3.] In Parts 1 and 2, the sequence $\varphi_0, \varphi_1, \dots, \varphi_{\alpha+1}$ is called a {\em derivation} (or {\em proof}) {\em of $\varphi$} and a {\em derivation of $\varphi$ from $\Gamma$}, respectively. Further, the length of this sequence is called the {\em length} of that derivation.
\end{itemize}
\end{definition}

As usual, we say that a formula $\varphi$ is {\em consistent} in $\mathcal{H}_{\mathsf{DPL}}$ if $\nvdash \neg\varphi$. Likewise, a set of formulas $w$ is  {\em consistent} if $w\nvdash\perp$.  Further, $w$ is {\em maximal consistent} if it is consistent and any set of formulas properly containing $w$ is inconsistent.  Finally, $w$ is {\em finitely consistent} whenever each finite subset of $w$ is consistent.

\begin{lemma} \label{prop}
	Let $\varphi$ and $\psi$ be formulas of $\mathcal{L}_{\mathsf{DPL}}$, and $\Gamma$ be a set of formula.
	\begin{itemize}
		\item[1.] $\vdash\!\bigcirc^n (\varphi \to\psi) \leftrightarrow (\bigcirc^n\varphi \to \bigcirc^n\psi)$, for all $n\in \mathbb{N}.\;\;\;\;\;\;\;\;\;\;\;$ (distribution)
		\item[2.]  $\vdash {\bigcirc}^n L_{r_1\dots r_k r} \varphi \to {\bigcirc}^n L_{r_1\dots r_k s} \varphi$, for all $n\in \mathbb{N}$, $r_1, \dots,  r_k$ and $s \leq r$.
	    \item[3.] If $\vdash\varphi \to\psi$ then $ \vdash L_r\varphi \to L_r\psi$, for all $r\in \mathbb{Q} \cap [0,1]$.
	    \item[4.] $\text{If } \Gamma \vdash\varphi \text { then }  \bigcirc\Gamma \vdash\bigcirc\varphi.\;\;\;\;\;\;\;\;\;\;\;\;\;\;\;\;\;$  (strong necessitation)
	\end{itemize}
\end{lemma}
\begin{proof}
	We shall only prove Part 4, as the rest are straightforward.
	\begin{itemize}
		\item[4.] Suppose that $\Gamma \vdash \varphi$. Then, by transfinite induction on the length of the derivation of $\varphi$ from $\Gamma$,  we show that $\bigcirc\Gamma \vdash \bigcirc\varphi$. The induction step in either of cases when $\varphi \in \Gamma$,  $\varphi$ is a theorem $\mathcal{H}_{\mathsf{DPL}}$, or $\varphi$ is obtained by an application of the rule MP, is easy. Therefore, we consider the case where $\varphi\equiv \theta \to{\bigcirc}^n L_{r_1\dots r_k r} \sigma$ is obtained from $\Gamma$ by an application of the rule GArch$_{\bigcirc^n, r}$. In this case, for all $s<r$, we have $\Gamma\vdash \theta \to{\bigcirc}^n L_{r_1\dots r_k s} \sigma$.
		By the induction hypothesis, $\bigcirc\Gamma\vdash \bigcirc (\theta \to{\bigcirc}^n L_{r_1\dots r_k s} \sigma)$, for all $s<r$. By Part 1, we obtain  $\bigcirc\Gamma\vdash\bigcirc\theta\to{\bigcirc}^{n+1} L_{r_1\dots r_k s} \sigma$, for all $s<r$. Using GArch$_{\bigcirc^n, r}$, we have $\bigcirc\Gamma\vdash\bigcirc\theta\to{\bigcirc}^{n+1} L_{r_1\dots r_k r} \sigma$.  Again by Part 1, it follows that $\bigcirc\Gamma\vdash\bigcirc(\theta\to{\bigcirc}^{n} L_{r_1\dots r_k r} \sigma)$, as required.
	\end{itemize}
	\vspace*{-4mm}
\end{proof}

\begin{theo} [Strong soundness of $\mathcal{H}_{\mathsf{DPL}}$] \label{sound1}
	The proof system $\mathcal{H}_{\mathsf{DPL}}$ is strongly sound for the class $\mathcal{DMP}$ of all dynamic Markov processes, i.e. if $\Gamma  \vdash\varphi$ then $\Gamma\vDash\varphi$.
\end{theo}
\begin{proof}
It suffices to show that each axiom of $\mathcal{H}_{\mathsf{DPL}}$ is valid in $\mathcal{DMP}$ and each rule preserves validity. Furthermore, by transfinite induction on the length of a derivation, we have to show if $\Gamma  \vdash\varphi$ then $\Gamma\vDash\varphi$. To this end, one should check that the rules MP and GArch$_{\bigcirc^n, r}$ preserve under the satisfiability relation. For example, to verify the claim for
GArch$_{\bigcirc^n, r}$, we should show that if $\mathfrak{M} = \langle{\Omega, \mathcal{A}, T, f, v}\rangle$ is a model and $w\in \Omega$ such that $\mathfrak{M}, w \vDash \Gamma$ and $\mathfrak{M}, w \vDash \theta\to{\bigcirc}^n L_{r_1\dots r_k s} \sigma$, for all $s<r$, then, $\mathfrak{M}, w \vDash \Gamma$ and $\mathfrak{M}, w \vDash \theta\to {\bigcirc}^n L_{r_1\dots r_k r} \sigma$. To see this, let $\mathfrak{M}, w \vDash \Gamma$ and $\mathfrak{M}, w \vDash \theta$. Then, by the induction hypothesis, we have $\mathfrak{M}, w \vDash {\bigcirc}^n L_{r_1\dots r_k s} \sigma$, for all $s<r$.
 Using the satisfiability relation,  $\mathfrak{M}, f^n(w) \vDash L_{r_1\dots r_k s}\sigma$, for all $s<r$. Hence by Fact \ref{continuity}, it follows that $\mathfrak{M}, f^n(w) \vDash L_{r_1\dots r_k r}\sigma$. Consequently, we obtain that $\mathfrak{M}, w \vDash {\bigcirc}^n L_{r_1\dots r_k r} \sigma$.
\end{proof}

\begin{theo} [Deduction theorem for $\mathsf{DPL}$] \label{ded1}
	Let $\varphi$ and $\psi$ be formulas of $\mathcal{L}_\mathsf{DPL}$ and $\Gamma$ a set of formulas. Then, $\Gamma, \varphi\vdash\psi$ iff $\Gamma\vdash\varphi\to \psi$.
\end{theo}
\begin{proof}
	The direction of right to left is obvious. For the other direction, we use transfinite induction on the length of the derivation of $\psi$ from $\Gamma\cup \{\varphi\}$. The induction step in either cases when $\psi =\varphi$, $\psi \in \Gamma$, $\psi$ is a theorem of $\mathcal{H}_{\mathsf{DPL}}$ or $\psi$ is obtained by an application of the rule MP, is standard. Thus, we consider the case where $\psi\equiv \theta \to{\bigcirc}^n L_{r_1\dots r_k r} \sigma$ is obtained from $\Gamma\cup \{\varphi\}$ by an application of  the rule GArch$_{\bigcirc^n, r}$. In this case,  for all $s< r$, we have $ \Gamma, \varphi\vdash\theta\to{\bigcirc}^n L_{r_1\dots r_k s} \sigma $.  By the induction hypothesis,  for all $s< r$, we obtain $\Gamma\vdash\varphi\to (\theta\to{\bigcirc}^n L_{r_1\dots r_k s} \sigma)$. So by the propositional validity, $\Gamma\vdash \varphi\wedge\theta\to{\bigcirc}^n L_{r_1\dots r_k s} \sigma$,  for all $s<r$. From this and using GArch$_{\bigcirc^n, r}$, it follows that $\Gamma\vdash\varphi\wedge\theta\to{\bigcirc}^n L_{r_1\dots r_k r} \sigma$.  Again by the propositional validity, we obtain $\Gamma\vdash\varphi\to(\theta\to{\bigcirc}^n L_{r_1\dots r_k r} \sigma)$, as desired.
\end{proof}
\subsection{Strong completeness of $\mathcal{H}_{\mathsf{DPL}}$} \label{Strong completeness of DPL}
In this subsection, we demonstrate that $\mathcal{H}_{\mathsf{DPL}}$ is strongly complete with respect to the class of all dynamic Markov processes based on standard Borel spaces. To this end, we will introduce special finitely consistent subsets of formulas, called {\em saturated sets}. Subsequently, we define a canonical model construction based on these sets, leading to strong completeness. Here, our approach is motivated by the methods, employed by Kozen {\em et al.} in \cite[Section 4]{Kozen2} to prove the strong completeness for Markovian logics. In the following, we denote the set of all formulas of $\mathcal{L}_\mathsf{DPL}$ by $\mathbb{F}$.

\begin{definition} \label{satu}
A set  $w \subseteq \mathbb{F}$ of formulas $\mathcal{L}_\mathsf{DPL}$ is  called {\em saturated} if
\begin{itemize}
\item[1.] $w$ is finitely consistent,
\item[2.] $w$ is negation complete,  i.e. for every formula $\varphi$ of $\mathcal{L}_\mathsf{DPL}$,  either $\varphi \in w$ or  $\neg\varphi \in w$,  and
\item[3.] $w$ has the Archimedean property, that is, for every formula $\varphi$ of $\mathcal{L}_\mathsf{DPL}$, $n\in \mathbb{N}$ and  $r_1, \dots, r_k, r\in \mathbb{Q} \cap [0,1]$, if $\{{\bigcirc}^n L_{r_1\dots r_k s} \varphi\;|\; s<r\}\subseteq w$ then ${\bigcirc}^n L_{r_1\dots r_k r} \varphi \in w$.
\end{itemize}
\end{definition}
Recall that a set of formulas is said to be  {\em maximal finitely consistent} if it satisfies Conditions 1 and 2  of the above definition.

\begin{pro} Let $w$ be a  saturated set of formulas of $\mathcal{L}_\mathsf{DPL}$. Then
 \begin{itemize}
\item[1.]  $w$ is closed under deduction, i.e. for every formula $\varphi$ of $\mathcal{L}_\mathsf{DPL}$, if $w\vdash\varphi$ then $\varphi \in w$;
\item[2.] $w$ is consistent.
\end{itemize}
\end{pro}
\begin{proof}
Part 1 can be proved by the transfinite induction on the length of the derivation of $\varphi$ from $w$. The only interesting case where the derivation ends with an application of the rule GArch$_{\bigcirc^n, r}$. In this case, the third condition of Definition \ref{satu} guarantees if $\varphi\equiv  \theta \to{\bigcirc}^n L_{r_1\dots r_k r} \sigma$ and $w\vdash\varphi$, then $ \theta\to{\bigcirc}^n L_{r_1\dots r_k r} \sigma \in w$.

Part 2 immediately follows from Part 1 and this fact that  $\perp\notin w$, since $w$ is finitely consistent.
\end{proof}

\begin{coro}
A set of formulas of $\mathcal{L}_\mathsf{DPL}$ is saturated if and only if it is maximal consistent.
\end{coro}
\begin{proof}
The direction of right to left is simple, since every maximal consistent set is closed under deduction. The other direction, immediately follows from the above proposition and the negation completeness property of the saturated sets.
\end{proof}

The following lemma is the key component of the proof of strong completeness theorem.

\begin{lemma} [Lindenbaum lemma] \label{lin}
Let $w$ be a consistent set of formulas of $\mathcal{L}_\mathsf{DPL}$. Then there exists a  saturated (or maximal consistent) set $w^*$  such that $w \subseteq w^*$.
\end{lemma}
\begin{proof}
Let  $\varphi_0, \; \varphi_1,\;   \varphi_2, \;\dots$ and $s_0, \; s_1,\;   s_2, \;\dots$ be two enumerations of all formulas of $\mathcal{L}_\mathsf{DPL}$ and all rational numbers in $[0, 1]$, respectively.  We inductively define a sequence
$\Gamma_0  \subseteq \Gamma_1  \subseteq  \dots \subseteq\Gamma_k  \subseteq \dots$ as follows:

 Set $\Gamma_0 := w$ and suppose that $\Gamma_k$ is already defined. Then, we put
\begin{equation*}
\Gamma_{k+1} := \left\{
\begin{array}{rl}
\Gamma_k \cup \{\varphi_k\}\;\;\;\;\; \;\;\;\;\;\;\;\;\;\;\;\;\; \;\;\;\;\;\;\;\;\;\;\;\;\;
& \text{if } \Gamma_k \vdash\varphi_k,\\
\Gamma_k \cup \{\neg\varphi_k\}\;\;\;\;\;\;\;\;\;\;\;\;\;\;\;\; \;\;\;\;\;\;\;\;\;\;\;\;\;\!
& \text{if }  \Gamma_k \nvdash\varphi_k \;\text{and $\varphi_k$ is not of}\\
& \text{the form}\;{\bigcirc}^n L_{r_1\dots r_k r} \theta,\\
\Gamma_k \cup \{\neg\varphi_k, \neg {\bigcirc}^n L_{r_1\dots r_k s_l} \theta\}\;\;\;
& \text{if }  \Gamma_k \nvdash\varphi_k \;\text{and}\;\varphi_k \equiv {\bigcirc}^n L_{r_1\dots r_k r} \theta\\
& \text{and $l$ is the least number such}\\
& \text{that $s_l< r$ and }\;   \Gamma_k \nvdash {\bigcirc}^n L_{r_1\dots r_k s_l} \theta.
\end{array} \right.
\end{equation*}
where $n\in \mathbb{N}$ and  $r_1, \dots, r_k, r, s_l\in \mathbb{Q} \cap [0,1]$. Notice that such $l$  in the third case always exists. As otherwise if  $\Gamma_k\vdash {\bigcirc}^n L_{r_1\dots r_k s_l} \theta$ for all $s_l<r$, then,  the rule GArch$_{\bigcirc^n, r}$ implies that $\Gamma_k\vdash {\bigcirc}^n L_{r_1\dots r_k r} \theta$, which is a contradiction.

Now, if we consider
\begin{equation*}
w^* := \bigcup_{k \in \omega}\Gamma_k
\end{equation*}
\\
then it is easy to show that $w^*$ satisfies all the desired properties.

 First, by  induction on $k$, it can be shown that all $\Gamma_k$ are consistent.  All cases are straightforward. The only interesting case is  when $\Gamma_{k+1} := \Gamma_k \cup \{\neg\varphi_k, \neg {\bigcirc}^n L_{r_1\dots r_k s_l} \theta\}$  where $\varphi_k \equiv {\bigcirc}^n L_{r_1\dots r_k r} \theta$ and $s_l<r$. Suppose that $\Gamma_{k+1}\vdash\bot$. Then by using the deduction theorem, $\Gamma_{k}\vdash {\bigcirc}^n L_{r_1\dots r_k r} \theta \vee {\bigcirc}^n L_{r_1\dots r_k s_l} \theta$. Now since by Part 2 of Lemma \ref{prop}, we have $\vdash{\bigcirc}^n L_{r_1\dots r_k r} \theta \to {\bigcirc}^n L_{r_1\dots r_k s_l} \theta$,
it follows that $\Gamma_{k}\vdash{\bigcirc}^n L_{r_1\dots r_k s_l} \theta$, which is a contradiction.

Clearly, by the above construction, $w^*$ is negation complete and has the Archimedean property. Therefore, it suffices to show that it is finitely consistent. To see this, let
$w'\subseteq w^*$ be finite and $w'\vdash\perp$. Then, we have $\vdash \neg \bigwedge w'$. Hence, for some $m\in \omega$, $ \neg \bigwedge w'\in \Gamma_m$. On the other hand, $w'\subseteq w^*$ implies that $w' \subseteq\Gamma_{m'}$ for some $m'\in \omega$. So it follows that  $\Gamma_k$ is inconsistent for $k= \text{max}\{m, m'\}$,  a contradiction.
\end{proof}

\begin{lemma}  \label{saturated}
For every pair of formulas $\varphi$ and $\psi$, we have $\vdash \varphi\to \psi$  if and only if for every  saturated set $w$, $\varphi \in w$ implies $\psi \in w$.
\end{lemma}
\begin{proof}
The direction of left to right is obvious, since $\varphi\to \psi\in w$. For the other direction, assume that  $\nvdash\varphi\to \psi$. Then $\varphi \land\neg \psi$ is consistent, since $\neg(\varphi\to \psi)\vdash\varphi \land\neg \psi$. So by the Lindenbaum lemma, there exists a saturated set $w^*$ such that $\varphi \land\neg \psi\in w^*$.  This implies that $\varphi \in w^*$ and $\psi\notin w^*$.
\end{proof}

Below, we show that the family of saturated sets forms a topological space.

\begin{definition} \label{topo}
	Suppose that $\Omega_c$ is the set  of  all  saturated sets of formulas of $\mathcal{L}_\mathsf{DPL}$. Then, let $\tau_c$ be the topology generated by the set $\mathcal{B}_c = \{[\varphi] \;|\;\varphi \in \mathbb{F}\}$ on $\Omega_c$ where for every formula $\varphi\in \mathbb{F}$, $[\varphi] :=  \{ w\in \Omega_c \;|\; \varphi\in w\}$. In this situation, the pair  $\langle{\Omega_c,  {\tau}_c}\rangle$ is called the {\em canonical topological space}.
\end{definition}

\begin{re} \label{subspace}
	One can easily check that $\mathcal{B}_c$, introduced in the above definition, forms a basis for $\tau_c$. Let $\Omega^*$ be the set of all maximal finitely consistent sets of formulas of $\mathcal{L}_\mathsf{DPL}$ and $\tau^*$ be the topology generated by the set $\mathcal{B}^* = \{[\varphi]^* \;|\;\varphi \in \mathbb{F}\}$ on $\Omega^*$ where for every formula $\varphi \in \mathbb{F}$, $[\varphi]^* :=  \{ w\in \Omega^* \;|\; \varphi\in w\}$. By the classical Stone representation theorem \cite{Stone}, the pair $\langle{\Omega^*,  {\tau}^*\!}\rangle$ is a Stone space, i.e., a compact totally disconnected Hausdorff space. Further, it is second countable since $\mathbb{F}$ and hence $\mathcal{B}^*$ are countable. Consequently, by Fact \ref{locally}, $\langle{\Omega^*,  {\tau}^*\!}\rangle$ is a Polish space. Note that for every formula $\varphi\in\mathbb{F}$, $[\varphi]= \Omega_c \cap [\varphi]^*$, hence the pair $\langle{\Omega_c,  {\tau}_c}\rangle$ is a subspace of  $\langle{\Omega^*,  {\tau}^*\!}\rangle$. For more details, see \cite[Section 6]{Kozen2013}.
\end{re}

Indeed, the following lemma shows that $\langle{\Omega_c,  {\tau}_c}\rangle$ also forms a Polish space. For an algebraic form of this lemma, see \cite[Proposition 13]{Kozen2013}.

\begin{lemma} \label{Polish}
The pair $\langle{\Omega_c,  {\tau}_c}\rangle$ is a Polish space.
\end{lemma}
\begin{proof}
By Remark \ref{subspace} and Fact \ref{G}, it suffices to show that $\Omega_c$ is a $G_\delta$ subset of $\Omega^*$.
To this end, for each $n\in\mathbb{N}$, consider
\[F_n= \{{\bigcirc}^n L_{r_1\dots r_k r} \varphi \;|\; \varphi\in \mathbb{F}\text{ and } r_1,\dots, r_k, r\in \mathbb{Q} \cap [0,1]\}.\]
Note that in $F_n$, the formulas of the form ${\bigcirc}^n L_{r_1\dots r_k r} \varphi$ can be parametrized by $r$. Hence, denote ${\bigcirc}^n L_{r_1\dots r_k r} \varphi \in F_n$ by $\psi^r$.
Further, the countability of $\mathbb{F}$ implies the countability of $F_n$, for each $n\in\mathbb{N}$. Now, let
\begin{align*}
O_{\psi^{r}} & = \{w^* \in \Omega^*\;|\; \psi^{r}\in w^* \text{ or } \exists s<r\; \;\psi^{s} \notin w^*\}\\
& = [\psi^{r}]^* \cup \bigcup_{s<r} [\neg \psi^{s}]^*.
\end{align*}
Then, by the Archimedean property of the saturated sets in $\Omega_c$, we have
\[\Omega_c = \bigcap_{n\in \mathbb{N}}  \bigcap_{\psi^{r}\in F_n} O_{\psi^{r}}.\]
It is easy to see that for each $\psi^{r}\in F_n$, $O_{\psi^{r}}$ is an open set in ${\tau}^*\!$. Consequently, $\Omega_c$ is a $G_\delta$ subset of $\Omega^*$.
\end{proof}

\begin{coro} \label{A}
Let  $\mathcal{A}_c$ be the Borel $\sigma$-algebra generated by $\tau_c$. Then, the pair $\langle{\Omega_c,  \mathcal{A}_c}\rangle$ is a standard Borel space.
\end{coro}
\begin{proof}
It is immediate by Lemma  \ref{Polish}.
\end{proof}

For a proof of  the following fact, see Lemmas 8 to 10 and Theorem 11 in \cite{Kozen2}.

\begin{fact}\label{T}
	\begin{itemize}
		\item[1.] It is easy to check that $\Omega_c \setminus [\varphi]= [\neg\varphi]$ and $[\varphi] \cap [\psi]=[\varphi \land \psi]$, for every $\varphi, \psi \in \mathbb{F}$. Therefore, the set $\mathcal{B}_c$ is an algebra on $\Omega_c$.
		\item[2.] For each $w\in \Omega_c$, the set function $\mu_c(w): \mathcal{B}_c\to [0, 1]$, defined as \[\mu_c(w)([\varphi]):= \sup{\{r\in \mathbb{Q} \cap [0,1] \;|\; L_r\varphi\in w\}} = \inf{\{r\in \mathbb{Q} \cap [0,1] \;|\; \neg L_r\varphi \in w\}},\] is finitely additive and continuous from above at the empty set.
		\item[3.] For each $w\in \Omega_c$,  $\mu_c(w)$ has an unique extension to a countably additive measure $T_c(w): \mathcal{A}_c\to [0, 1]$. Further, it is a probability measure.
		\item[4.] $T_c: \Omega_c\times \mathcal{A}_c\to [0, 1]$ is a measurable function.
	\end{itemize}
\end{fact}

\begin{lemma} \label{f}
Consider $f_c: \Omega_c\to 2^{\mathbb{F}}$ defined as ${f_c}(w): = \{ \varphi\in \mathbb{F} \;|\;\bigcirc\varphi \in w\}$. Then,
\begin{itemize}
\item[1.] For each $w \in \Omega_c$, we have $f_c(w) \in \Omega_c$, and hence $f_c$ can be considered as a function from $\Omega_c$ to $\Omega_c$.
\item[2.] $f_c$ is a measurable function.
\end{itemize}
\end{lemma}
\begin{proof}
For Part 1, we have
\begin{itemize}
\item[1.] $f_c(w)$ is finitely consistent. Suppose not. Then, there is a finite subset  $w'$ of $f_c(w)$ such that $w' \vdash \perp$. Then, $\vdash \neg \bigwedge w'$.
By the rule Nec$_\bigcirc$, we have $\vdash \bigcirc\neg \bigwedge w'$. So by the axioms Func and Conj$_\bigcirc$, we obtain that $\vdash \neg \bigwedge \bigcirc w'$. So $\bigcirc w' \vdash \perp$. But, this contradicts the finite consistency of  $w$, since  $\bigcirc w'$ is a finite subset of $w$.
\item[2.] $f_c(w)$ is negation complete. Assume that $\varphi\notin f_c(w)$. Then, we have $\bigcirc\varphi\notin w$. So $\neg\bigcirc\varphi\in w$, since $w$ is negation complete. By the axiom Func$_\bigcirc$,  $\bigcirc\neg\varphi\in w$. This implies that $\neg\varphi\in f_c(w)$. Similarly, it can be shown that if  $\neg\varphi\in f_c(w)$ then $\varphi\notin f_c(w)$.
\item[3.] $f_c(w)$ has the Archimedean property. Suppose that   $\{{\bigcirc}^n L_{r_1\dots r_k s} \varphi\;|\; s<r\}\subseteq f_c(w)$. Then,   $\{{\bigcirc}^{n+1} L_{r_1\dots r_k s} \varphi\;|\; s<r\}\subseteq w$ by the definition of $f_c$.  Now since $w$ has the Archimedean property,  so ${\bigcirc}^{n+1} L_{r_1\dots r_k r} \varphi \in w$. Consequently, ${\bigcirc}^{n} L_{r_1\dots r_k r} \varphi \in f_c(w)$.
\end{itemize}

For Part 2,  it suffices to show that for each  $\varphi\in\mathbb{F}$,  $f_{c}^{-1}([\varphi]) = [\bigcirc\varphi]$.
Let $w\in f_{c}^{-1}([\varphi]) $. Then, we have $f_{c}(w)\in [\varphi]$. Hence, $\varphi\in  f_{c}(w)$. This yields that $\bigcirc \varphi\in w$. So, we have  $w\in [\bigcirc\varphi]$.
The other direction can be shown similarly.
\end{proof}

In the view of  Corollary \ref{A}, Fact \ref{T} and Lemma \ref{f}, we introduced a dynamic Markov model which can be exploited as a tool to establish strong completeness theorem.

\begin{definition}[Canonical model for $\mathsf{DPL}$] \label{cpmodel1}
	The \textnormal{canonical model} for $\mathsf{DPL}$ is the tuple $\mathfrak{M}_c = \langle{\Omega_c, \mathcal{A}_c, T_c, f_c, v_c}\rangle$ where
	\begin{itemize}
		\item[1.]  $\langle{\Omega_c, \mathcal{A}_c, T_c}\rangle$ is the Markov process obtained from Fact \ref{T};
		\item[2.] $f_c: \Omega_c \to \Omega_c$ defined as ${f_c}(w) = \{ \varphi\in \mathbb{F} \;|\;\bigcirc\varphi \in w\}$ for each $w\in \Omega_c$;
		\item[3.] $v_c$ is the valuation defined as $v_c(p) = \{ w\in \Omega_c \;|\;p\in w\}$ for each $p \in \mathbb{P}$.
	\end{itemize}
\end{definition}

\begin{lemma}[Truth lemma] \label{truth2}
 Consider the canonical model  ${\mathfrak{M}}_c$. Then for every formula $\varphi$ of $\mathcal{L}_\mathsf{DPL}$ and $w\in\Omega_c$,  we have  ${\mathfrak{M}}_c,  w \vDash\varphi$ iff $\varphi \in w.$
\end{lemma}
\begin{proof}
Both direction can be simultaneously proved by induction on complexity of a formula $\varphi$. The proofs are straightforward and follow from the satisfiability relation of each logical complexities of formulas.
\end{proof}

The truth lemma yields the strong completeness theorem for  $\mathcal{H}_{\mathsf{DPL}}$.

\begin{theo}[Strong completeness of $\mathcal{H}_{\mathsf{DPL}}$] \label{scompletef}
$\mathcal{H}_{\mathsf{DPL}}$ is strongly complete with respect to the class $\mathcal{MS}$ of all dynamic Markov processes based on standard Borel spaces.
\end{theo}
\begin{proof}
Let $\Gamma\cup\{\varphi\}$ be a set of formulas of $\mathcal{L}_\mathsf{DPL}$ such that $\Gamma\nvdash\varphi$. Then $\Gamma\cup\{\neg\varphi\}$ is consistent. Therefore, by Lindenbaum lemma, there is a  saturated set $w^*$  such that $ \Gamma\cup\{\neg\varphi\}\subseteq w^*$.  But by truth lemma, we have ${\mathfrak{M}}_c, w^*\vDash\gamma$ for each $\gamma\in \Gamma$ and ${\mathfrak{M}}_c, w^*\nvDash\varphi$. So  we have a dynamic Markov model based on a standard Borel space with a node in which  $\Gamma$ is satisfied and  $\varphi$ is not satisfied, that is, $\Gamma\;{\nvDash}_{\mathcal{MS}}\;\varphi$.
\end{proof}
\subsection{The logic of some special classes of dynamic Markov processes} \label{special}
In this subsection, we examine the logics of some important subclasses of dynamic Markov processes, introduced in Definition \ref{frames}.

\begin{definition}  \label{Axiom}
Consider the following axiom schemes and inference rule:
	\begin{align*}
	&\textnormal{M} && L_r \bigcirc\varphi \leftrightarrow \bigcirc L_r\varphi &&\;\;\;\;\;\;\;\;\;\;\;\;\;\;\;\;\;\;\;\;\;\;\;\;\;\;\;\!\;\;\text{(measure-preserving)}\;\;\;\;\;\;\;\;\;\;\;\;\;\;\;\;\\
	&\textnormal{Id} && \bigcirc\varphi \leftrightarrow \varphi &&\;\;\;\;\;\;\;\;\;\;\;\;\;\;\;\;\;\;\;\;\;\;\;\;\;\;\;\;\;\text{(identity)}\;\;\;\;\;\;\;\;\;\;\;\;\;\;\;\;\\
	&\textnormal{H$_1$}&& L_r \varphi \to L_1 L_r\varphi  &&\;\;\;\;\;\;\;\;\;\;\;\;\;\;\;\;\;\;\;\;\;\;\;\;\;\;\;\;\;\text{(Harsanyi1)}\\
	&\textnormal{H$_2$}&& \neg L_r \varphi \to L_1\neg L_r\varphi &&\;\;\;\;\;\;\;\;\;\;\;\;\;\;\;\;\;\;\;\;\;\;\;\;\;\;\;\;\;\text{(Harsanyi2)}\\
	&\textnormal{GArch$_r$} && \frac{\displaystyle{\;\{\psi\to L_{r_1\dots r_k s}\varphi \;|\; s< r\}\;}}{\displaystyle{\psi \to L_{r_1\dots r_k r} \varphi}}&& \text{\;\;\;\;\;\;\;\;\;\;\;\;\;\;\;\;\;\;\;\;\;\;\;\;\;}\;\;\;\;\;\;\;\;\;\;\;\;\;\;\;\;\;\;\;\;\;\;\;\;\;\;\;\;\;\;\;\;\;\;\;\;\;\;
	\end{align*}
	Then,
	\begin{itemize}
		\item[1.] let $\mathcal{H}_{\mathsf{M}}$ be the proof system whose axioms consist of all axioms of $\mathcal{H}_{\mathsf{DPL}}$ together with axiom scheme M, and the rules MP, Nec$_{L_1}$, Nec$_\bigcirc$ and GArch$_r$.
		\item[2.] Likewise, consider the proof system $\mathcal{H}_{\mathsf{Pure}} := \mathcal{H}_{\mathsf{M}} + \text{Id}$.
		\item [3.] Finally, let  $\mathcal{H}_{\mathsf{ADS}}:= \mathcal{H}_{\mathsf{M}} + \text{H}_1 + \text{H}_2$
	\end{itemize}
\end{definition}

\begin{re}
	Because of presence of axiom M in the above logics, it is easy to check that  the scheme 	$ L_{r_1\dots r_k r} {\bigcirc}^n \varphi \leftrightarrow {\bigcirc}^n L_{r_1\dots r_k r}\varphi$ is a theorem of these logics, for each $n \in \mathbb{N}$ and $r_1, \dots, r_k, r \in \mathbb{Q}\cup [0, 1]$. Hence, the rule  GArch$_{\bigcirc^n, r}$ can be derived in these logics.
\end{re}

Below,  we show that the proof systems $\mathcal{H}_{\mathsf{M}}$, $\mathcal{H}_{\mathsf{Pure}}$ and $\mathcal{H}_{\mathsf{ADS}}$ are strongly sound and complete for the class $\mathcal{M}$ of all measure-preserving dynamic Markov processes, the class $\mathcal{P}{ure}$ of all purely probabilistic dynamic Markov processes and the class $\mathcal{ADS}$ of all abstract dynamical systems,  respectively.  The proof of strong soundness for  all mentioned logics as well as strong completeness of  $\mathcal{H}_{\mathsf{M}}$ and $\mathcal{H}_{\mathsf{Pure}}$ are straightforward. But, the proof for strong completeness of $\mathcal{H}_{\mathsf{ADS}}$ is more elaborate.

\begin{theo} \label{MP and Pure}
	\begin{itemize}
		\item[1.] $\mathcal{H}_{\mathsf{M}}$ is strongly sound and complete for the class $\mathcal{M}$.
		\item[2.] $\mathcal{H}_{\mathsf{Pure}}$ is strongly sound and complete for the class $\mathcal{P}{ure}$.
	\end{itemize}
\end{theo}
\begin{proof}
The proof of strong soundness in both cases are left for the reader to verify.  For the proof of strong completeness, notice that in the canonical model $\mathfrak{M}_c= \langle{\Omega_c, \mathcal{A}_c, T_c, f_c, v_c}\rangle$, introduced in Definition \ref{cpmodel1}, the function $f_c$ is measure-preserving (resp. identity) if axiom M (resp. Id) is valid in $\mathfrak{M}_c$.
\end{proof}

\begin{re} \label{P}
On the basis of the second part of the above theorem,   $\mathcal{H}_{\mathsf{Pure}}$ is essentially nothing but {\em probability logic}; see e.g. \cite{Kozen2, Zhou2009}.
\end{re}

\begin{theo} \label{PS}
	 $\mathcal{H}_{\mathsf{ADS}}$ is strongly sound and complete for the class $\mathcal{ADS}$.
\end{theo}
\begin{proof}
We only proof the strong completeness. Let $\Gamma\cup\{\varphi\}$ be a set of formulas of  $\mathcal{L}_\mathsf{DPL}$ and $\Gamma \nvdash \varphi$ in $\mathcal{H}_{\mathsf{ADS}}$. So $\Gamma \nvdash \varphi$ in $\mathcal{H}_\mathsf{DPL}$. Hence, in the light of Theorem \ref{scompletef}, we have $\Gamma \nvDash \varphi$. Now, if we take the canonical model  $\mathfrak{M}_c= \langle{\Omega_c, \mathcal{A}_c, T_c, f_c, v_c}\rangle$, then there exists $w_0\in \Omega_c$ such that while $\mathfrak{M}_c, w_0\vDash \Gamma$,  we have $\mathfrak{M}_c, w_0\nvDash\varphi$.

By  \cite[Proposition 1]{Meier}, since the axioms H$_1$ and H$_2$ are valid in the  canonical model $\mathfrak{M}_c$,  this model enjoys the {\em Harsanyi property} H, namely,
	\begin{itemize}
		\item[1.] $\{v\in \Omega_c\;|\;T_c(v, [\varphi])\neq T_c(w, [\varphi])\} \in \mathcal{A}_c$, and
		\item[2.] $T_c(w) (\{v\in \Omega_c\;|\;T_c(v, [\varphi])\neq T_c(w, [\varphi])\}) =0$
	\end{itemize}
 for each $w\in \Omega_c$  and $\varphi \in \mathbb{F}$. But, since $\mathbb{F}$ is countable, it follows that for each $w\in \Omega_c$,
 	\begin{itemize}
 	\item[1.] $\bigcup_{\varphi \in \mathbb{F}}  \{v\in \Omega_c\;|\;T_c(v, [\varphi])\neq T_c(w, [\varphi])\} \in \mathcal{A}_c$, and
 	\item[2.] $T_c(w) (\bigcup_{\varphi \in \mathbb{F}}  \{v\in \Omega_c\;|\;T_c(v, [\varphi])\neq T_c(w, [\varphi])\}) =0$.
 \end{itemize}
Now set
\[A_0^{w_0}= \bigcup_{\varphi \in \mathbb{F}}  \{v\in \Omega\;|\;T_c(v, [\varphi])\neq T_c(w_0, [\varphi])\}\]
and  for each $n\in \mathbb{N}$, inductively define $A_{n+1}^{w_0}:= f^{-1}_c (A_n^{w_0})$. Then since $f_c $ is measure-preserving, this implies that $T_c(w_0)(A_n^{w_0}) =0 $. Hence, if we take $A_\infty^{w_0}:= \bigcup_{n\in \mathbb{N}}  A_n^{w_0}$, then $T_c(A_\infty^{w_0})= 0$ and $T_c(w_0)(\Omega_c \setminus A_\infty^{w_0}) = 1$.  Let   $\Omega_c^{'} = \Omega_c \setminus A_\infty^{w_0}$. Then, it is easy to check that $w_0 \in \Omega_c^{'}$ and $f_c(\Omega_c^{'}) \subseteq \Omega_c^{'}$. Moreover, for each $w\in \Omega_c^{'}$, we have that $T_c(w_0, [\varphi]) = T_c(w, [\varphi])$. But, since both $T_c (w_0, .)$ and $T_c(w, .)$ are probability measures on $\mathcal{A}_c$, and $\mathcal{A}_c$ is generated by the  countable $\mathcal{B}_c= \{[\varphi]\;|\;\varphi\in \mathbb{F}\}$, we have that $T_c (w_0, A) = T_c (w, A)$ for each $A\in \mathcal{A}_c$. Put $\mathcal{A}_c^{'}:= \{ A\cap \Omega_c^{'}\;|\;  A\in \mathcal{A}_c\}$. Note that $\mathcal{A}_c^{'}$ is a $\sigma$-algebra and $\mathcal{A}_c^{'}\subseteq  \mathcal{A}_c$. Therefore, if we put
$ T_c^{'}(w, .) = T_c(w_0, .) \restriction_{\mathcal{A}_c^{'}}$, then $T_c^{'}$ is probability measure and that $T_c^{'}(w_0, .)= T_c^{'}(w, .)$ for each $w\in \Omega_c^{'}$. Now consider $\mathfrak{M}_c^{'}= \langle{\Omega_c^{'}, \mathcal{A}_c^{'}, T_c^{'}, f_c^{'}, v_c^{'}}\rangle$ where $f_c^{'}:= f_c \restriction_{\Omega_c^{'}}$ and $v_c^{'}(p):= v_c(p)\cap \Omega_c^{'}$ for each propositional variable  $p\in \mathbb{P}$. It is clear that the underlying dynamic Markov process of $\mathfrak{M}_c^{'}$ is a dynamic probability space.

Now, by induction on the complexity of formulas in $\mathcal{L}_\mathsf{DPL}$, we prove that
\[\mathfrak{M}_c, w\vDash \theta \;\text{ iff }\;\mathfrak{M}_c^{'}, w\vDash \theta\]
	for each $w \in \Omega_c^{'}$. For example, to show the induction step for $\theta= L_r\psi$, suppose that $\mathfrak{M}_c, w\vDash L_r\psi$. Then, $T_c(w, [\![\psi]\!]_{\mathfrak{M}_c})\geq r$. But, by induction hypothesis for $\psi$, we have $[\![\psi]\!]_{\mathfrak{M}^{'}_c} = [\![\psi]\!]_{\mathfrak{M}_c} \cap \Omega_c^{'}$. Furthermore, $ T_c(w, [\![\psi]\!]_{\mathfrak{M}_c})= T_c(w, [\![\psi]\!]_{\mathfrak{M}_c} \cap \Omega_c^{'}) = T_c^{'}(w, [\![\psi]\!]_{\mathfrak{M}_c^{'}})$. Thus, $\mathfrak{M}_c^{'}, w\vDash L_r\psi$.
The other direction can be proved in a similar way.
	
Now, once the above claim	is established, we have that  $\mathfrak{M}_c^{'}, w_0\vDash \Gamma$ and $\mathfrak{M}_c^{'}, w_0\nvDash\varphi$, and the proof is complete.
\end{proof}

\begin{re}
Notice that the canonical models in the proof of the above strong completeness, are based on some standard Borel spaces.
\end{re}

\section{Infinitary dynamic probability logic} \label{IDPL}
In this section, we introduce the infinitary dynamic probability logic by adding countable conjunctions and disjunctions to the language $\mathcal{L}_\mathsf{DPL}$. We present a weakly complete Hilbert-style axiomatization for this logic and show that when the logic is restricted to its countable fragments, the proposed axiomatization is strongly complete with respect to the class $\mathcal{DMP}$ of all dynamic Markov processes.
\subsection{Syntax and semantics} \label{Syntax and semantics 2}
As in the previous section, we assume that $\mathbb{P}$ is a countable set of propositional variables. The language $\mathcal{L}_{\mathsf{DPL}_{\omega_1}}\!$ of {\em infinitary dynamic probability logic} extends the language $\mathcal{L}_{\mathsf{DPL}}$
by adding infinitary conjunction $\bigwedge_{i\in I} \varphi _i$ where $\{\varphi _i\;|\;i\in I\}$ is a countable set of formulas. So the language $\mathcal{L}_{\mathsf{DPL}_{\omega_1}}\!$  is recursively defined by the following grammar:
\begin{align*}
\varphi ::&=\ \   p\mid \neg\varphi \mid\bigwedge_{i\in I} \varphi _i\mid L_r \varphi \mid \bigcirc\varphi
\end{align*}
where $p \in \mathbb{P}$,  $\{\varphi _i\;|\;i\in I\}$ is a countable set of formulas, and $r\in \mathbb{Q} \cap [0,1]$. We define the usual abbreviations for the propositional connectives $\top$, $\perp$, $\to$, $\leftrightarrow$, $\bigvee$ and the probability operators $M_r$ and $L_{r_1\dots r_k}$. Notice that when $I=\{1, 2\}$, then $\bigwedge_{i\in I} \varphi _i= \varphi_1 \wedge \varphi_2$ and hence $\bigvee_{i\in I} \varphi _i = \varphi_1 \vee \varphi_2$. So obviously the set $\mathbb{F}$ of all formulas of
$\mathcal{L}_{\mathsf{DPL}}$ is a subset of the set of all formulas of $\mathcal{L}_{\mathsf{DPL}_{\omega_1}}\!$. Accordingly, a formula $\varphi$ in $\mathbb{F}$ is called {\em finitary formula} while the other formulas in $\mathcal{L}_{\mathsf{DPL}_{\omega_1}}\!$ are called {\em infinitary formulas}.
Moreover, it is not hard to check that the cardinality of the set of all formulas of $\mathcal{L}_{\mathsf{DPL}_{\omega_1}}$ is $2^{\aleph_0}$.

Given a dynamic Markov model $ \mathfrak{M} = \langle{\Omega, \mathcal{A}, T, f, v}\rangle$, the satisfiability relation for the new formulas of $\mathcal{L}_{\mathsf{DPL}_{\omega_1}}\!$ is defined inductively as follows:
\[\mathfrak{M}, w \vDash \bigwedge_{i\in I}\varphi_i \;\text{ iff }\; \mathfrak{M}, w \vDash \varphi_i \text{ for all } i\in I. \;\;\;\;\;\;\;\;\;\;\;\;\;\;\;\;\;\;\;\;\;\;\;\;\;\;\;\;\;\;\;\;\;\;\;\;\;\;\;\;\;\;\;\;\;\;\;\;\;\;\;\;\;\;\;\;\;\;\;\;\;\;\;\;\;\;\]

\begin{definition}
	By {\em infinitary dynamic probability logic}, or $\mathsf{DPL}_{\omega_1}$, we mean the set of all valid formulas of $\mathcal{L}_{\mathsf{DPL}_{\omega_1}}$ over  the class $\mathcal{DMP}$ of all  dynamic Markov processes.
\end{definition}
\subsection{The axiomatization of $\mathsf{DPL}_{\omega_1}$} \label{The axiomatization of IDPL}
The following Hilbert-style axiomatization, denoted by $\mathcal{H}_{\mathsf{DPL}_{\omega_1}}$, is given for $\mathsf{DPL}_{\omega_1}$. We subsequently show that $\mathcal{H}_{\mathsf{DPL}_{\omega_1}}$  is strongly sound, Theorem \ref{sound2}, and weakly complete, Corollary \ref{com1}, with respect to the class of all dynamic Markov processes.

\begin{definition}
	The proof system $\mathcal{H}_{\mathsf{DPL}_{\omega_1}}$ has the following axiom schemes:
	\begin{align*}
	& \textnormal{ITaut}&&\text{All  infinitary propositional tautologies }\\
	&\textnormal{FA$_1$}&& L_0 \perp &&\text{(finite additivity)}\\
	&\textnormal{FA$_2$}&&L_r \neg\varphi \to \neg L_s\varphi, \;\text{if}\;\; r+s > 1 &&\text{}\\
	&\textnormal{FA$_3$}&& L_r (\varphi\land \psi) \land L_s(\varphi\land \neg\psi) \to L_{r+s}\varphi, \;\text{if}\;\; r+s \leq 1&&\text{}\\
	&\textnormal{FA$_4$}&& \neg L_r (\varphi\land \psi) \land \neg L_s(\varphi\land \neg\psi) \to \neg L_{r+s}\varphi, \;\text{if}\;\; r+s \leq 1&&\text{}\\
	&\textnormal{Mono}&&   L_1 (\varphi\to\psi) \to (L_r\varphi\to L_r \psi)  &&\text{(monotonicity)}\\
	&\textnormal{Arch}&& \bigwedge_{s< r} L_s \varphi \to L_r \varphi &&\text{(Archimedean)}\\
	&\textnormal{Cont}&& \bigwedge_{0<n\in \mathbb{N}}\; \bigvee_{k\in \mathbb{N}}\; \neg L_{\frac{1}{n}} (\bigwedge_{i\leq k} \varphi_i \land \neg \bigwedge_{i \in \mathbb{N}} \varphi_i) &&\text{(continuity)}\\
	&\textnormal{Func$_\bigcirc$}\;&&\bigcirc\neg\varphi\leftrightarrow \neg\bigcirc\varphi&&\text{(functionality)}\\
	&\textnormal{Conj$_\bigcirc$}\;&&\bigcirc\bigwedge_{i \in I}\varphi_i \leftrightarrow \bigwedge_{i \in I}\bigcirc\varphi_i, \text{ if } |I|<\aleph_1 &&\text{(conjunction)}
	\end{align*}
	and the inference rules:
	\begin{align*}
	&\textnormal{MP} && \frac{\displaystyle{\;\psi\to\varphi\;\;\;\psi\;}}{\displaystyle{\varphi}}\;\;\;\;\;\;\;\;\;\;\;\;\;\;\;\;\;\;\;\;\;\;\;\;\;\;\;\;\;\;\;\;\;\;\;\;\;\;\;\;\;\;\;\;\;\;\;\;\;\;\;\;\;\;\;\;\;\;\;\;\;\;\;\;\;\;&& \text{(modus ponens)} \;\;\\
	&\textnormal{Conj} && \frac{\displaystyle{\;\{\psi\to \varphi_i\; |\;i \in I\}\;}}{\displaystyle{\psi\to \bigwedge_{i\in I} \varphi _i}}, \text{ for } |I| < \aleph_1 && \text{(conjunction)}\\
&\textnormal{Arch$_{\bigcirc^n, r}$}&&\frac{\displaystyle{\;\{\psi\to{\bigcirc}^n L_ s \varphi \;|\; s< r\}\;}}{\displaystyle{\psi \to {\bigcirc}^n L_r \varphi}}, \text{ for } \; n\in \mathbb{N} \;\;\;\;\;\;\;\;\;\;\;\;\;\;\;\;\;\;\;\;\;\;\;\; && \text{(Archimedean)}\\
	&\textnormal{Nec$_{L_1}$} && \frac{\displaystyle{\;\varphi\;}}{\displaystyle{L_1 \varphi}} && \text{(necessitation)}\\
	&\textnormal{Nec$_\bigcirc$} && \frac{\displaystyle{{\;\varphi\;}}}{\displaystyle{\bigcirc\varphi}} && \text{}
	\end{align*}
\end{definition}

Note that here all formulas are in $\mathcal{L}_{\mathsf{DPL}_{\omega_1}}$. The {\em continuity axiom} syntactically expresses the property of the continuity from above at the empty set for a finite measure. For more details,  we refer the reader to \cite[Remark 3]{Bara}.

\begin{definition}\label{Theorem2}
\begin{itemize}
\item[1.] Similar to $\mathsf{DPL}$, one can define a {\em theorem} $\varphi$ of $\mathcal{H}_{\mathsf{DPL}_{\omega_1}}$ (denoted by $\vdash_{\omega_1}\varphi$), by assigning a sequence $\varphi_0, \varphi_1, \dots, \varphi_{\alpha+1}$ ($\alpha+1 < \omega_1$) of the formulas of $\mathcal{L}_{\mathsf{DPL}_{\omega_1}}$ with $\varphi= \varphi_{\alpha+1}$ and for each  $\beta \leq{\alpha+1}$, either $\varphi_\beta$ is an axiom or is derived by applying one of the inference rules on some preceding formulas of the sequence.
\item[2.] Likewise, a formula $\varphi$ is {\em derivable from a set $\Gamma$ of assumptions} in $\mathcal{H}_{\mathsf{DPL}_{\omega_1}}$ (denoted by $\Gamma  \vdash_{\omega_1}\!\varphi$), if there exists a sequence $\varphi_0, \varphi_1, \dots, \varphi_{\alpha+1}$ ($\alpha+1 < \omega_1$) of the formulas of $\mathcal{L}_{\mathsf{DPL}_{\omega_1}}$  with $\varphi = \varphi_{\alpha+1}$ and for each $\beta \leq{\alpha+1}$, either $\varphi_\beta$ is a member of $\Gamma$, or is a theorem of  $\mathcal{H}_{\mathsf{DPL}_{\omega_1}}$, or is derived by applying one of the inference rules, other than  Nec$_{L_1}$ and Nec$_\bigcirc$, on some preceding formulas of the sequence.
\end{itemize}
\end{definition}

The notions of a {\em derivation} (or {\em proof}) and {\em length of a derivation} as well as {\em consistency} of a set of formulas can be defined similar to $\mathsf{DPL}$. Further, the strong soundness theorem can be proved similar to Theorem \ref{sound1}.

\begin{theo} [Strong soundness of $\mathcal{H}_{\mathsf{DPL}_{\omega_1}}$] \label{sound2}
	Let $\Gamma\cup \{\varphi\}$ be a set of formulas of $\mathcal{L}_{\mathsf{DPL}_{\omega_1}}$. Then, $\Gamma \vdash_{\omega_1}\!\varphi$ implies $\Gamma\vDash\varphi$. 
\end{theo}
\begin{proof}
	The proof runs similar to the proof of Theorem \ref{sound1}.
\end{proof}

\begin{theo} [Deduction theorem for $\mathsf{DPL}_{\omega_1}$] \label{ded}
	Let $\varphi$ and $\psi$ be formulas of $\mathcal{L}_{\mathsf{DPL}_{\omega_1}}$ and $\Gamma$ a set of formulas. Then, $\Gamma, \varphi\vdash_{\omega_1}\!\psi$ iff $\Gamma \vdash_{\omega_1}\!\varphi\to \psi$.
\end{theo}
\begin{proof}
	The proof is similar to the proof of Theorem \ref{ded1}. The only new part is the case where $\psi\equiv\theta\to \bigwedge_{i\in I}\sigma_i$ for $|I|< \aleph_1$,  is obtained from $\Gamma\cup\{\varphi\}$ by an application of the rule Conj. In this case, $\Gamma, \varphi \vdash_{\omega_1}\! \theta\to\sigma_i$ for all $i\in I$. Using the induction hypothesis, $\Gamma  \vdash_{\omega_1}\! \varphi\to (\theta\to \sigma_i)$ for all $i \in I$. This implies that $\Gamma  \vdash_{\omega_1}\! \varphi\wedge\theta\to \sigma_i$ for all $i \in I$. By the rule Conj, it follows that $\Gamma \vdash_{\omega_1}\!\varphi\wedge\theta\to \bigwedge_{i\in I}\sigma_i$.  So we have $\Gamma \vdash_{\omega_1}\!\varphi\to(\theta\to \bigwedge_{i\in I}\sigma_i)$.
\end{proof}

In the next lemmas, we give some of the main properties of the derivation relation $\vdash_{\omega_1}$ and some examples of theorems of $\mathcal{H}_{\mathsf{DPL}_{\omega_1}}\!$ which will be used later.

\begin{lemma} \label{prop1}
	Let $\varphi$ and $\psi$ be formulas of  $\mathcal{L}_{\mathsf{DPL}_{\omega_1}}$ and  $\Gamma$ and $\Delta$ be set of formulas. Then,
	\begin{itemize}
        \item[1.]  $\Gamma \cup \Delta\vdash_{\omega_1}\!\varphi$ iff $\Gamma, \bigwedge \Delta \vdash_{\omega_1}\! \varphi$.
        \item[2.] If $\Gamma \vdash_{\omega_1}\!\delta$ for all $\delta\in \Delta$ and $\Delta \vdash_{\omega_1}\!\varphi$, then $\Gamma\vdash_{\omega_1}\!\varphi$.
       \item[3.] If $\Gamma, \psi \vdash_{\omega_1}\! \varphi$ and $\Gamma, \neg\psi \vdash_{\omega_1}\! \varphi$, then $\Gamma \vdash_{\omega_1}\! \varphi$.
        \item[4.] If $\Gamma \vdash_{\omega_1}\! \varphi$, then $\{\psi \to \gamma\;|\; \gamma\in\Gamma \} \vdash_{\omega_1}\! \psi\to\varphi$.
      \item[5.] If $\Gamma \vdash_{\omega_1}\!\varphi$ then $\bigcirc\Gamma \vdash_{\omega_1}\!\bigcirc\varphi$.
	    \item[6.] If $\vdash_{\omega_1}\!\varphi \to\psi$ then $ \vdash_{\omega_1}\! L_r\varphi \to L_r\psi$, for all $r\in \mathbb{Q} \cap [0,1]$.
	    \item[7.] $\{\varphi_i\to \psi\; |\;i \in I\} \vdash_{\omega_1}\! \bigvee_{i\in I} \varphi_i \to \psi$, where $\{\varphi_i\; |\;i \in I\}$ is a countable set of formulas.
	     \item[8.] $\{\psi\to L_{\frac{1}{n}} (L_s \varphi\land \neg L_r \varphi) \; | \;s<r\} \vdash_{\omega_1}\!\neg \psi $, for $0< n\in \mathbb{N}$.
	     \item[9.] $\{\psi\to\bigcirc^n  L_{r_1 \dots r_k s} \varphi\; | \;s<r\} \vdash_{\omega_1}\! \psi \to \bigcirc^n L_{r_1 \dots r_k r}  \varphi$, for $n\in \mathbb{N}$.	
	\end{itemize}
\end{lemma}
\begin{proof}
	We shall prove Parts 8 and 9. The proof of the other parts is straightforward.
	 \begin{itemize}
	 \item[8.]  By the axiom Cont, we have that
	     \[\vdash_{\omega_1}\! \bigwedge_{0<n\in \mathbb{N}}\;\bigvee_{0< k\in \mathbb{N}}\; \neg L_{\frac{1}{n}}(\bigwedge_{0< i \leq k}  L_{r\dot{-} \frac{1}{i}} \varphi \wedge \neg  \bigwedge_{0< i \in \mathbb{N}}  L_{r\dot{-} \frac{1}{i}} \varphi).\footnote{Here $r\dot{-}\frac{1}{i}$ is truncated minus, that is, $r\dot{-}\frac{1}{i}=$ if $r\leq\frac{1}{i}$ then $0$ else $r-\frac{1}{i}$.}\]
On the other hand, by the axiom Arch and this fact that $\vdash_{\omega_1}\! L_r \varphi \to L_s \varphi$ for all $s< r$, we obtain $\vdash_{\omega_1}\!  L_r \varphi \leftrightarrow  \bigwedge_{0< i\in \mathbb{N}}  L_{r\dot{-} \frac{1}{i}} \varphi$ and $\vdash_{\omega_1}\!   L_{r\dot{-} \frac{1}{k}} \varphi \leftrightarrow  \bigwedge_{0< i \leq k}  L_{r\dot{-} \frac{1}{i}} \varphi$, for all $0< k \in \mathbb{N}$.  From this, we can conclude that
	  \[\vdash_{\omega_1}\! \bigwedge_{0<n\in \mathbb{N}}\;\bigvee_{0< k\in \mathbb{N}}\; \neg L_{\frac{1}{n}}( L_{r\dot{-} \frac{1}{k}} \varphi \wedge \neg L_r \varphi).\]
	Therefore,  for all $0< n \in \mathbb{N}$, we have  $\vdash_{\omega_1}\! \neg \bigwedge_{0< k\in \mathbb{N}}\; L_{\frac{1}{n}}( L_{r\dot{-} \frac{1}{k}} \varphi \wedge \neg L_r \varphi)$. From this and using Part 4, it is easy to see that
\[\{\psi\to L_{\frac{1}{n}} (L_s \varphi\land \neg L_r \varphi) \; | \;s<r\} \vdash_{\omega_1}\!\neg \psi.\]
	 \item[9.] By Part 8, it can be shown that $\{\psi\to L_t L_s \varphi\; | \;s<r\} \vdash_{\omega_1}\! \psi \to L_t L_r \varphi$, for all $t, s, r\in \mathbb{Q} \cap [0,1]$. For a proof, see \cite[Theorem 4.3]{Iko}. Using this and induction on $k$,  it is not hard to see that  $\{ L_{r_1 \dots r_k s} \varphi\; | \;s<r\} \vdash_{\omega_1}\!  L_{r_1 \dots r_k r}  \varphi$, for all $r_1 \dots r_k\in \mathbb{Q}\cap [0,1]$.  By Part  5,  we can derive that $\{ \bigcirc^n L_{r_1 \dots r_k s} \varphi\; | \;s<r\} \vdash_{\omega_1}\! \bigcirc^n L_{r_1 \dots r_k r}  \varphi$,  for any $n\in \mathbb{N}$. So, by Part 4, we obtain
	 $\{\psi\to\bigcirc^n  L_{r_1 \dots r_k s} \varphi\; | \;s<r\} \vdash_{\omega_1}\! \psi \to \bigcirc^n L_{r_1 \dots r_k r}  \varphi$.
	 \end{itemize}
	 \vspace*{-4mm}
\end{proof}

\begin{lemma} \label{theo1}
	Let $\varphi$ and $\psi$ be formulas of  $\mathcal{L}_{\mathsf{DPL}_{\omega_1}}$ and $\{\varphi_i\;|\; i\in I\}$ be a countable set of formulas. Then, we have
	\begin{itemize}
        \item[1.]  $\vdash_{\omega_1}\!\bigcirc^n (\varphi \to\psi) \leftrightarrow (\bigcirc^n\varphi \to \bigcirc^n\psi)$, for all $n\in \mathbb{N}$.
		\item[2.] $\vdash_{\omega_1}\!\bigcirc\bigvee_{i \in I}\varphi_i \leftrightarrow \bigvee_{i \in I}\bigcirc\varphi_i.$
		\item[3.] $\vdash_{\omega_1}\! L_r (\bigwedge_{i \in \mathbb{N}}\varphi_i) \leftrightarrow \bigwedge_{k\in \mathbb{N}} L_r(\bigwedge_{i\leq k}\varphi_i).$
    	\item[4.] $\vdash_{\omega_1}\! L_r (\bigvee_{i \in \mathbb{N}}\varphi_i) \leftrightarrow \bigwedge_{0<n\in \mathbb{N}}\; \bigvee_{k\in \mathbb{N}} L_{r\dot{-}\frac{1}{n}}(\bigvee_{i\leq k}\varphi_i).$
	\end{itemize}
\end{lemma}
\begin{proof}
	We will prove Parts  3 and 4, the rest are straightforward.
	 \begin{itemize}
	 \item[3.]  For the left to right direction,  for all $k \in \mathbb{N}$, we have $\vdash_{\omega_1}\!\bigwedge_{i \in \mathbb{N}}\varphi_i \to \bigwedge_{i\leq k} \varphi_i$. From this and using Part 6 of Lemma \ref{prop1}, we obtain $\vdash_{\omega_1}\! L_r(\bigwedge_{i \in \mathbb{N}}\varphi_i) \to L_r(\bigwedge_{i\leq k} \varphi_i)$  for all $k \in \mathbb{N}$. So by the rule Conj, it follows that  $\vdash_{\omega_1}\! L_r (\bigwedge_{i \in \mathbb{N}}\varphi_i) \to \bigwedge_{k\in \mathbb{N}} L_r(\bigwedge_{i\leq k}\varphi_i).$
	
	For the other direction,  by the axiom Arch, we have
	  \[\vdash_{\omega_1}\! \neg L_r (\bigwedge_{i \in \mathbb{N}}\varphi_i) \to \bigvee_{0< n\in \mathbb{N}}  \neg L_{r \dot{-} \frac{1}{n}}(\bigwedge_{i \in \mathbb{N}}\varphi_i).\]
	  Further, by the axiom Cont,   $\vdash_{\omega_1}\! \bigwedge_{0<n\in \mathbb{N}}\;\bigvee_{k\in \mathbb{N}}  \neg L_{\frac{1}{n}} (\bigwedge_{i\leq k} \varphi_i \land \neg \bigwedge_{i \in \mathbb{N}} \varphi_i)$.  On the other hand,  by the axiom FA$_4$, for all  $0<n\in \mathbb{N}$ and $k \in \mathbb{N}$, we have
	 \[\vdash_{\omega_1}\! \neg L_{r \dot{-} \frac{1}{n}}(\bigwedge_{i \in \mathbb{N}}\varphi_i) \land  \neg L_{\frac{1}{n}} (\bigwedge_{i\leq k} \varphi_i \land \neg \bigwedge_{i \in \mathbb{N}} \varphi_i) \to \neg L_r (\bigwedge_{i\leq k} \varphi_i).\]
	By this fact that $\vdash_{\omega_1}\!\neg L_r (\bigwedge_{i\leq k} \varphi_i) \to \bigvee_{k\in  \mathbb{N}}  \neg L_r (\bigwedge_{i\leq k} \varphi_i)$, we can conclude that  for all  $0<n\in \mathbb{N}$ and $k \in \mathbb{N}$,
	 \[\vdash_{\omega_1}\! \neg L_{r \dot{-} \frac{1}{n}}(\bigwedge_{i \in \mathbb{N}}\varphi_i) \land  \neg L_{\frac{1}{n}} (\bigwedge_{i\leq k} \varphi_i \land \neg \bigwedge_{i \in \mathbb{N}} \varphi_i) \to \bigvee_{k\in  \mathbb{N}}  \neg L_r (\bigwedge_{i\leq k} \varphi_i).\]
	 By Part 7 of Lemma \ref{prop1} and the fact $\vdash_{\omega_1}\!  \bigvee_{k\in \mathbb{N}} (\psi \land \theta_k) \leftrightarrow (\psi \land \bigvee_{k\in \mathbb{N}} \theta_k) $, we obtain
	  \[\vdash_{\omega_1}\! \neg L_{r \dot{-} \frac{1}{n}}(\bigwedge_{i \in \mathbb{N}}\varphi_i) \land  \bigvee_{k\in  \mathbb{N}}  \neg L_{\frac{1}{n}} (\bigwedge_{i\leq k} \varphi_i \land \neg \bigwedge_{i \in \mathbb{N}} \varphi_i) \to \bigvee_{k\in  \mathbb{N}}  \neg L_r (\bigwedge_{i\leq k} \varphi_i),\]
	  for all  $0<n\in \mathbb{N}$. Using the fact $\vdash_{\omega_1}\! \bigwedge_{n\in \mathbb{N}} \sigma_n \to \sigma_n$ for all $n\in \mathbb{N}$,  it is easy to see  that
	   \[\vdash_{\omega_1}\! \neg L_{r \dot{-} \frac{1}{n}}(\bigwedge_{i \in \mathbb{N}}\varphi_i) \land  \bigwedge_{0< n\in \mathbb{N}} \;\bigvee_{k\in  \mathbb{N}}  \neg L_{\frac{1}{n}} (\bigwedge_{i\leq k} \varphi_i \land \neg \bigwedge_{i \in \mathbb{N}} \varphi_i) \to \bigvee_{k\in  \mathbb{N}}  \neg L_r (\bigwedge_{i\leq k} \varphi_i),\]
	  for all  $0<n\in \mathbb{N}$. Again by Part 7 of Lemma \ref{prop1}, it can be shown that
	   \[\vdash_{\omega_1}\! \bigvee_{0< n\in \mathbb{N}} \neg L_{r \dot{-} \frac{1}{n}}(\bigwedge_{i \in \mathbb{N}}\varphi_i) \land \bigwedge_{0< n\in \mathbb{N}} \; \bigvee_{k\in  \mathbb{N}}  \neg L_{\frac{1}{n}} (\bigwedge_{i\leq k} \varphi_i \land \neg \bigwedge_{i \in \mathbb{N}} \varphi_i) \to \bigvee_{k\in  \mathbb{N}}  \neg L_r (\bigwedge_{i\leq k} \varphi_i).\]
	   Consequently, we have $ \vdash_{\omega_1}\!\neg L_r (\bigwedge_{i \in \mathbb{N}}\varphi_i) \to \bigvee_{k\in  \mathbb{N}}  \neg L_r (\bigwedge_{i\leq k} \varphi_i)$. This means that
	   \[\vdash_{\omega_1}\! \bigwedge_{k\in  \mathbb{N}}  L_r (\bigwedge_{i\leq k} \varphi_i)    \to L_r (\bigwedge_{i \in \mathbb{N}}\varphi_i).\]
	 \item[4.] Using this fact that $\vdash_{\omega_1}\!  \neg L_r \psi \to L_{1-r} \neg \psi$, we obtain
	\[\vdash_{\omega_1}\!  \neg \bigwedge_{0<n\in \mathbb{N}}\;  \bigvee_{k \in \mathbb{N}}\; L_{r\dot{-} \frac{1}{n}} (\bigvee_{i\leq k} \varphi_i) \to  \bigvee_{0<n\in \mathbb{N}}\; \bigwedge_{k \in \mathbb{N}}\; L_{1-r+\frac{1}{n}} (\bigwedge_{i\leq k} \neg \varphi_i).\]
	From this and Part 3, it follows that
	\[\vdash_{\omega_1}\!  \neg \bigwedge_{0<n\in \mathbb{N}}\;  \bigvee_{k \in \mathbb{N}}\; L_{r\dot{-} \frac{1}{n}} (\bigvee_{i\leq k} \varphi_i) \to  \bigvee_{0<n\in \mathbb{N}}\; L_{1-r+\frac{1}{n}} (\bigwedge_{i  \in \mathbb{N}} \neg \varphi_i).\]
	 Now, by the axiom FA$_2$, we can conclude that
	 \[\vdash_{\omega_1}\!  \neg \bigwedge_{0<n\in \mathbb{N}}\;  \bigvee_{k \in \mathbb{N}}\; L_{r\dot{-} \frac{1}{n}} (\bigvee_{i\leq k} \varphi_i) \to \neg \bigwedge_{0<n\in \mathbb{N}}\; L_{r\dot{-}\frac{1}{2n}} (\bigvee_{i  \in \mathbb{N}}  \varphi_i).\]
Thus, by this fact that $\vdash_{\omega_1}\!  L_r \psi \to \bigwedge_{s<r} L_s \psi$, we get 	
  \[\vdash_{\omega_1}\!  \neg \bigwedge_{0<n\in \mathbb{N}}\;  \bigvee_{k \in \mathbb{N}}\; L_{r\dot{-} \frac{1}{n}} (\bigvee_{i\leq k} \varphi_i) \to \neg  L_{r} (\bigvee_{i  \in \mathbb{N}}  \varphi_i).\]
	
	For the converse,  we have $\vdash_{\omega_1}\! L_{r \dot{-} \frac{1}{n}} (\bigvee_{i \leq k} \varphi_i) \to L_{r \dot{-} \frac{1}{n}} (\bigvee_{i \in \mathbb{N}} \varphi_i)$,  for all $k\in \mathbb{N}$ and $0<n \in \mathbb{N}$. From this and using Part 7 of Lemma \ref{prop1},  it follows that   for all $0<n \in \mathbb{N}$,
	\[\vdash_{\omega_1}\! \bigvee_{k\in  \mathbb{N}} L_{r \dot{-} \frac{1}{n}} (\bigvee_{i \leq k} \varphi_i) \to L_{r \dot{-} \frac{1}{n}} (\bigvee_{i \in \mathbb{N}} \varphi_i).\]
Further, it is easy to see that
	 $\vdash_{\omega_1}\! \bigwedge_{0<n\in \mathbb{N}}\; \bigvee_{k\in  \mathbb{N}} L_{r \dot{-} \frac{1}{n}} (\bigvee_{i \leq k} \varphi_i) \to L_{r \dot{-} \frac{1}{n}} (\bigvee_{i \in \mathbb{N}} \varphi_i)$  for all $0<n \in \mathbb{N}$. So, by the rule Conj, we obtain
	 \[\vdash_{\omega_1}\! \bigwedge_{0<n\in \mathbb{N}}\; \bigvee_{k\in  \mathbb{N}} L_{r \dot{-} \frac{1}{n}} (\bigvee_{i \leq k} \varphi_i) \to \bigwedge_{0<n\in \mathbb{N}} L_{r \dot{-} \frac{1}{n}} (\bigvee_{i \in \mathbb{N}} \varphi_i).\] Now by the axiom Arch, we can conclude that
	 \[\vdash_{\omega_1}\! \bigwedge_{0<n\in \mathbb{N}} \;\bigvee_{k\in  \mathbb{N}} L_{r \dot{-} \frac{1}{n}} (\bigvee_{i \leq k} \varphi_i) \to L_{r} (\bigvee_{i \in \mathbb{N}} \varphi_i).\]
	 \end{itemize}
	 \vspace*{-4mm}
\end{proof}

\begin{definition}
Two formulas $\varphi$ and $\psi$ are called {\em provably contradictory} if $\vdash_{\omega_1} \neg(\varphi \wedge \psi)$. Moreover, the formulas $\varphi_1, \dots, \varphi_n$ are said to be {\em pairwise provably  contradictory} if any two of which are provably contradictory.
\end{definition}

\begin{lemma} \label{prop2}
	Let two formulas $\varphi$ and $\psi$ be provably  contradictory and the formulas $\varphi_1, \dots, \varphi_{n+1}$ be pairwise provably  contradictory. Then, we have the following expressions:
	\begin{itemize}
		\item[1.] $\vdash_{\omega_1}\! L_r\varphi \wedge L_s\psi \to L_{r+s} (\varphi \vee \psi)$, for $r+s\leq 1$.
		\item[2.] $\vdash_{\omega_1}\! \neg L_r\varphi \wedge \neg L_s\psi \to \neg L_{r+s} (\varphi \vee \psi)$, for $r+s\leq 1$.
		\item[3.] $\vdash_{\omega_1}\! L_r(\varphi \vee \psi)\leftrightarrow \bigwedge_{s< r} \;(L_{s}\varphi \vee L_{r-s}\psi).$
\item[4.] $\vdash_{\omega_1}\!  L_r(\varphi_1 \vee\dots \vee \varphi_n\vee \varphi_{n+1})\leftrightarrow \bigwedge_{s_1+ \dots +s_{n}< r} \;(L_{s_1} \varphi_1 \vee \dots \vee L_{s_{n}} \varphi_{n} \vee L_{r-(s_1+ \dots +s_{n})}\varphi_{n+1}).$
	\end{itemize}
\end{lemma}
\begin{proof}
	We shall only prove Parts 3 and 4. For a proof of Parts 1 and 2, see  \cite[Lemma~3.6]{Zhou2009}.
	 \begin{itemize}
	 \item[3.] For the left to right direction, by Part 2  we have that $\vdash_{\omega_1}\! L_r(\varphi \vee \psi)\to L_{s}\varphi \vee L_{r-s}\psi$, for all $s<r$.  Using the rule Conj, it follows that
	$\vdash_{\omega_1}\! L_r(\varphi \vee \psi)\to \bigwedge_{s< r} \;(L_{s}\varphi \vee L_{r-s}\psi).$
	
	For the converse, it is enough to show that
	\begin{itemize}
	\item[(i)] $\neg L_0 \varphi \vee L_r \varphi,\; \bigwedge_{s< r} \;(L_{s}\varphi \vee L_{r-s}\psi) \vdash_{\omega_1}\! L_r(\varphi \vee \psi)$, and
	\item[(ii)] $L_0 \varphi \wedge \neg L_r \varphi, \; \bigwedge_{s< r} \;(L_{s}\varphi \vee L_{r-s}\psi) \vdash_{\omega_1}\! L_r(\varphi \vee \psi)$.
	\end{itemize}
	Since by (i) and (ii) and Part 3 of Lemma \ref{prop1}, we have $\bigwedge_{s< r} \;(L_{s}\varphi \vee L_{r-s}\psi) \vdash_{\omega_1}\! L_r(\varphi \vee \psi)$. So, $\vdash_{\omega_1}\!\bigwedge_{s< r} \;(L_{s}\varphi \vee L_{r-s}\psi)\to L_r(\varphi \vee \psi)$.

The proof of (i) is easy and follows from the fact that $\vdash_{\omega_1}\!   L_0 \varphi$, since in this case we have  $\neg L_0 \varphi \vee L_r \varphi\vdash_{\omega_1}\!   L_r (\varphi \vee \psi)$.
			
For the proof of (ii), first, by induction on $n$, we show that for all $0<n\in \mathbb{N}$,
\[L_0 \varphi \wedge \neg L_r \varphi \vdash_{\omega_1}\!  \bigvee_{k< 2^n} (L_{\frac{rk}{2^n}}\varphi \wedge \neg L_{\frac{r(k+1)}{2^n}} \varphi).\]
For the base case $n=1$, we have $L_0 \varphi \wedge \neg L_r \varphi \vdash_{\omega_1}\! (L_0 \varphi \wedge \neg L_r \varphi) \wedge (L_{\frac{r}{2}} \varphi \vee \neg L_{\frac{r}{2}})$. From this, it is easy to see that $L_0 \varphi \wedge \neg L_r \varphi \vdash_{\omega_1}\! (L_0 \varphi \wedge  \neg L_{\frac{r}{2}}) \vee (L_{\frac{r}{2}} \varphi \vee \neg L_r \varphi)$, as desired.  Now, suppose that the claim is true for some $n >1$, we show that the claim also holds for $n+1$. By the induction hypothesis,  we have
	\[L_0 \varphi \wedge \neg L_r \varphi \vdash_{\omega_1}\!  \bigvee_{2k< 2^{n+1}} (L_{\frac{r(2k)}{2^{n+1}}}\varphi \wedge \neg L_{\frac{r(2k+2)}{2^{n+1}}} \varphi).\]
	From this and using the fact $\vdash_{\omega_1}\! L_{\frac{r(2k+1)}{2^{n+1}}}\varphi \vee \neg L_{\frac{r(2k+1)}{2^{n+1}}} \varphi$, we obtain that
	\[L_0 \varphi \wedge \neg L_r \varphi \vdash_{\omega_1}\!  \bigvee_{2k< 2^{n+1}} \big(L_{\frac{r(2k)}{2^{n+1}}}\varphi \wedge \neg L_{\frac{r(2k+2)}{2^{n+1}}} \varphi \wedge (L_{\frac{r(2k+1)}{2^{n+1}}}\varphi \vee \neg L_{\frac{r(2k+1)}{2^{n+1}}} \varphi)\big).\]
This implies that 	
	\[L_0 \varphi \wedge \neg L_r \varphi \vdash_{\omega_1}\!  \bigvee_{2k< 2^{n+1}} \big(( L_{\frac{r(2k)}{2^{n+1}}}\varphi \wedge \neg L_{\frac{r(2k+1)}{2^{n+1}}} \varphi) \vee (L_{\frac{r(2k+1)}{2^{n+1}}}\varphi \wedge \neg L_{\frac{r(2k+2)}{2^{n+1}}} \varphi)\big).\]
	From this, it is not hard to see that
	\[L_0 \varphi \wedge \neg L_r \varphi \vdash_{\omega_1}\!  \bigvee_{k< 2^{n+1}} (L_{\frac{rk}{2^{n+1}}}\varphi \wedge \neg L_{\frac{r(k+1)}{2^{n+1}}} \varphi).\]
Now, Let $\theta = L_0 \varphi \wedge \neg L_r \varphi \wedge \bigwedge_{s< r} \;(L_{s}\varphi \vee L_{r-s}\psi)$. Then, for all $0<n\in \mathbb{N}$,
	 \[\theta \vdash_{\omega_1}\! \bigvee_{k< 2^n} (L_{\frac{rk}{2^n}}\varphi \wedge \neg L_{\frac{r(k+1)}{2^n}} \varphi).\]
	 Further, for all $s<r$, we have $\theta \vdash_{\omega_1}\! \neg L_s\varphi \to L_{r-s} \psi$. So if $s= \frac{r(k+1)}{2^n}$, then for all $0<n\in \mathbb{N}$, we have
	 \[\theta \vdash_{\omega_1}\! \neg L_{\frac{r(k+1)}{2^n}}\varphi \to L_{r-{\frac{r(k+1)}{2^n}}} \psi.\]
	 Therefore, for all $0<n\in \mathbb{N}$ and $k<2^n$,
	 \[\theta, (L_{\frac{rk}{2^n}}\varphi \wedge \neg L_{\frac{r(k+1)}{2^n}} \varphi)\vdash_{\omega_1}\! L_{\frac{rk}{2^n}}\varphi \wedge  L_{r-{\frac{r(k+1)}{2^n}}} \psi,\]
	 and by Part 1,
	  \[\theta, (L_{\frac{rk}{2^n}}\varphi \wedge \neg L_{\frac{r(k+1)}{2^n}} \varphi)\vdash_{\omega_1}\! L_{r-{\frac{r}{2^n}}}(\varphi \vee \psi).\]
	 Hence, by Part 7 of Lemma \ref{prop1},
	 \[\theta, \bigvee_{k< 2^n} (L_{\frac{rk}{2^n}}\varphi \wedge \neg L_{\frac{r(k+1)}{2^n}} \varphi)\vdash_{\omega_1}\! L_{r-{\frac{r}{2^n}}}(\varphi \vee \psi),\]
	and so
	 \[\theta\vdash_{\omega_1}\! L_{r-{\frac{r}{2^n}}}(\varphi \vee \psi).\]
	Using the rule Conj,  we have
	 \[\theta \vdash_{\omega_1}\! \bigwedge_{0<n\in \mathbb{N}}L_{r-{\frac{r}{2^n}}}(\varphi \vee \psi).\]
	  By the axiom Arch, this implies that $\theta \vdash_{\omega_1}\! L_r(\varphi\vee\psi)$.
	 \item[4.] This can be proved by Part 3 and  induction on $n\geq 2$.
	 \end{itemize}
	 \vspace*{-4mm}
\end{proof}
\subsection{Introducing $\Delta_\infty$}
In the following, we introduce a special class of formulas $\Delta_\infty = \bigcup_{\alpha<{ (2^{\aleph_0})^+} } \Delta_\alpha$. This class plays an important role in proving completeness theorems \ref{com2} and \ref{com1}.

\begin{definition}
Let $\nu = (2^{\aleph_0})^+$. We inductively define a sequence $(\Delta_\alpha)_{\alpha <\nu}$ of sets of formulas of $\mathcal{L}_{\mathsf{DPL}_{\omega_1}}\!$ as follows:
\begin{itemize}
\item Let $\Delta_0$ be the set of all formulas of $\mathcal{L}_\mathsf{DPL}$, i.e. $\Delta_0:=\mathbb{F}$;
\item $\Delta_\alpha := \bigcup_{\beta<\alpha}  \Delta_\beta$ for each non-zero limit ordinal $\alpha <\nu $;
\item Let  the set
\[\Gamma_\alpha :=\Delta_\alpha \cup\{\bigwedge_{i\in I} \varphi _i, \bigvee_{i\in I} \varphi _i\;|\;  \varphi _i \in  \Delta_\alpha \text{ and } |I|< \aleph_1\}.\]
Then $\Delta_{\alpha+1}$ is defined as the set of all finite Boolean combinations of the members of $\Gamma_\alpha$. Hence, if $\{\varphi_i\;|\;i\in I\}$ is a countable subset of $\Delta_\alpha$, then $\bigwedge_{i\in I} \varphi _i$ and $\bigvee_{i\in I} \varphi _i$ are in $\Delta_{\alpha+1}$. Further, if $\varphi$ and $\psi$ are in $\Delta_{\alpha+1}$, then so are $\neg \varphi$, $\varphi \land \psi$ and $\varphi \vee \psi$ in $\Delta_{\alpha+1}$.
\end{itemize}
Finally, put $\Delta_\infty := \bigcup_{\alpha<\nu } \Delta_\alpha$.
\end{definition}

\begin{re} \label{repre1}
\begin{itemize}
\item[1.] It can be easily seen by transfinite induction that $\Delta_\alpha$ is closed under finite Boolean combinations for each $\alpha<\nu$.
\item[2.]  Further, it is known that every finite Boolean combinations of members of $\Gamma_\alpha$ is provably equivalent to a finite disjunctive normal form of members of  $\Gamma_\alpha$.
Therefore, if $\varphi \in \Delta_{\alpha+1}$ then there exists a formula of the form $\bigvee_{k=1}^{m} \bigwedge_{l=1}^{n_k} \varphi_{l}^{(k)}$ in $\Delta_{\alpha+1}$ such that for each $1 \leq k \leq m$ and $1\leq  l \leq  n_k$,  either the formula $\varphi_{l}^{(k)}\in \Gamma_\alpha$ or  $\neg \varphi_{l}^{(k)}\in \Gamma_\alpha$ and
$\vdash_{\omega_1}\!\varphi \leftrightarrow \bigvee_{k=1}^{m} \bigwedge_{l=1}^{n_k} \varphi_{l}^{(k)}$. Hence, as $\Delta_\alpha$  is closed under finite Boolean combinations, each $\varphi_{l}^{(k)}$  is either of the form $\bigwedge_{p\in \mathbb{N}} \sigma_{p}^{(k, l)}$ or $\bigvee_{p\in \mathbb{N}} \sigma_p^{(k, l)}$ where $\sigma_p^{(k, l)}\in \Delta_\alpha$ for each  $1 \leq k \leq m$, $1\leq  l \leq  n_k$ and $p\in \mathbb{N}$.
\item[3.]  Also note that $\varphi_{l}^{(k)}$'s can be chosen in such way that they are pairwise provably contradictory.
\item[4.] Finally, by the infinite distribution law\footnote{That is $\vdash_{\omega_1}\! \varphi\wedge \bigvee_{i\in I} \psi_i \leftrightarrow \bigvee_{i\in I} (\varphi\wedge \psi_i)$.}, each $\varphi \in \Delta_{\alpha+1}$ is provably equivalent with a formula of the form $\bigvee_{k=1}^{m} \gamma_k$ such that  $\gamma_k= \bigvee_{i\in \mathbb{N}} \bigwedge_{j\in \mathbb{N}} \psi_{i,j}^{(k)}$ with $\psi_{i, j}^{(k)}\in \Delta_\alpha$ and $\gamma_k$'s are  pairwise provably  contradictory  for each  $1 \leq k \leq m$.
\end{itemize}
\end{re}

Next we show that the set $\Delta_\infty$ is provably closed under the probability and next modal operators, that is, for each formula $\varphi \in \Delta_\infty$ there exist  formulas $\theta, \psi_r  \in \Delta_\infty$ such that $\vdash_{\omega_1}\!\bigcirc\varphi \leftrightarrow \theta$ and $\vdash_{\omega_1}\! L_r\varphi \leftrightarrow \psi_r$ for each $r\in \mathbb{Q} \cap [0,1]$.

\begin{lemma} \label{closure properties}
	$\Delta_\infty$ is closed under the probability and next operators.
\end{lemma}
\begin{proof}
The statement is proved by transfinite induction on $\alpha$. Showing the statement for the $\bigcirc$ operator is easy, since this operator commutes with all Boolean operators including infinite conjunctions and disjunctions. So by representation of a formula $\varphi \in \Delta_{\alpha+1}$, in Part 4 of Remark \ref{repre1}, and induction hypotheses, one can show that  there exists a formula $\theta\in \Delta_\infty$ such that
$\vdash_{\omega_1}\!\bigcirc\varphi \leftrightarrow \theta$. The proof for the $L_r$ operator is more elaborate. So we present the proof in detail.

Let $\varphi \in \Delta_{\alpha+1}$ and assume that for each $\sigma\in \Delta_\alpha$ and each $s \in \mathbb{Q}\cap [0, 1]$, we have that the formula
$L_s \sigma$  is provably equivalent to a formula in $\Delta_\infty$. Then, by Part 4 of Remark \ref{repre1}, $\varphi$ is provably equivalent with a formula of the form $\bigvee_{k=1}^{m} \gamma_k$ such that  $\gamma_k= \bigvee_{i\in \mathbb{N}} \bigwedge_{j\in \mathbb{N}} \psi_{i,j}^{(k)}$ with $\psi_{i, j}^{(k)}\in \Delta_\alpha$ and $\gamma_k$'s are pairwise provably contradictory  for each  $1 \leq k \leq m$. Using Part 4 of Lemma \ref{prop2}, we have
\[\vdash_{\omega_1}\! L_r\varphi \leftrightarrow  \bigwedge_{s_1+ \dots +s_{m-1} < r} \; (L_{s_1}\gamma_1 \vee \dots \vee  L_{s_{m-1}} \gamma_{m-1} \vee L_{r-(s_1+ \dots +s_{m-1})} \gamma_m)\]
where $s_1, \dots, s_{m-1}\in \mathbb{Q} \cap [0,1]$. To summarize, let us denote  $r-(s_1+ \dots +s_{m-1})$ by $s_m$. So, it is enough to show the induction hypothesis for each $\gamma_k$. By  Part 4 of Lemma \ref{theo1},  for each $1 \leq  k \leq  m$, we have that
\[ \vdash_{\omega_1}\! L_{s_k}\gamma_k  \leftrightarrow \bigwedge_{0<n\in \mathbb{N}}\; \bigvee_{l\in \mathbb{N}}\;L_{s_{k} \dot{-} \frac{1}{n}}(\bigvee_{i\leq l}  \; \bigwedge_{j\in \mathbb{N}} \psi_{i,j}^{(k)} ).\]

Now by applying distributivity finite disjunctions with infinite conjunctions, we can rearrange formulas in $(\bigvee_{i\leq l}  \; \bigwedge_{j\in \mathbb{N}} \psi_{i,j}^{(k)})$ to find a formula of the form
$ \bigwedge_{p\in \mathbb{N}}  \theta_{p}$ such that  $\theta_p\in \Delta_\alpha$ and
\[\vdash_{\omega_1}\! (\bigvee_{i\leq l}  \; \bigwedge_{j\in \mathbb{N}} \psi_{i,j}^{(k)})\leftrightarrow\bigwedge_{p\in \mathbb{N}}  \theta_{p}.  \]
Therefore, it is enough to verify the induction hypothesis for infinite conjunctions of member of $\Delta_\alpha$. But by Part 3 of Lemma \ref{theo1}, we have that
	\[ \vdash_{\omega_1}\!  L_t (\bigwedge_{p\in \mathbb{N}} \theta_p) \leftrightarrow  \bigwedge_{q\in \mathbb{N}} L_t (\bigwedge_{p\leq q} \theta_p)  .\]
	Now since  $ (\bigwedge_{p\leq q} \theta_p)\in \Delta_\alpha$,  by the  induction hypothesis, we have  $L_t  (\bigwedge_{p\leq q} \theta_p) \in \Delta_\infty$. Consequently, the induction hypothesis holds for infinite conjunctions members of $\Delta_\alpha$, and then the proof is complete.
\end{proof}

\begin{lemma} \label{equi}
	For each formula $\varphi$ of $\mathcal{L}_{\mathsf{DPL}_{\omega_1}}$, there is $\varphi' \in \Delta_\infty$ such that $\vdash_{\omega_1}\!\varphi\leftrightarrow\varphi' $.
\end{lemma}
\begin{proof}
	The proof immediately follows from induction on the complexity of $\varphi$ and the preceding lemma.
\end{proof}

The following lemma shows that any saturated subset of $\mathbb{F}$ (or $\Delta_0$) decides about every formulas in $\Delta_\infty$.

\begin{lemma} \label{negcompleteF}
	For any saturated  set $w$  of formulas in $\mathbb{F}$ and $\varphi \in  \Delta_\infty$,  we have the following properties:
	\begin{itemize}
		\item[1.] Either $ w \vdash_{\omega_1}\!\varphi $ or $w\vdash_{\omega_1}\! \neg\varphi$, but not both.
		\item[2.] $w\vdash_{\omega_1}\!\varphi$ iff $w\cup \{\varphi\}$ is consistent in $\mathcal{H}_{\mathsf{DPL}_{\omega_1}}$.
	\item[3.] For any $\varphi \in \mathbb{F}$,  $w\vdash_{\omega_1}\! \varphi$ iff $\varphi\in w$.
	\end{itemize}
\end{lemma}
\begin{proof}
	\begin{itemize}
		\item[1.]  As before, we use transfinite induction to show the statement holds for each $\varphi \in \Delta_\alpha$ with $\alpha<\nu$.
		
		The proofs for  $\alpha = 0$ as well as for a non-zero limit ordinal  $\alpha$ are immediate.
		
		Now suppose that the claim is true for $\alpha<\nu$ and let $\varphi \in \Delta_{\alpha+ 1}$. Then by Part 4 of Remark \ref{repre1}, $\varphi$ is provably equivalent with a formula of the form $\bigvee_{k=1}^{m} \gamma_k$ such that for each  $1 \leq k \leq m$, $\gamma_k= \bigvee_{i\in \mathbb{N}} \bigwedge_{j\in \mathbb{N}} \psi_{i,j}^{(k)}$ with $\psi_{i, j}^{(k)}\in \Delta_\alpha$. If for some $1 \leq k \leq m$, $w \vdash_{\omega_1}\!  \gamma_k$, then  $w \vdash_{\omega_1}\! \varphi$. Otherwise, suppose $w \nvdash_{\omega_1}\!  \gamma_k$ for each $1 \leq k \leq m$. Then by Part 7 of Lemma \ref{prop1} and the rule Conj, for each $1 \leq k \leq m$ and $i\in \mathbb{N}$, there exists $j\in \mathbb{N}$ such that $w \nvdash_{\omega_1}\!   \psi_{i, j}^{(k)}$ and equivalently, $w\vdash_{\omega_1}\!  \neg\psi_{i, j}^{(k)}$ by the induction hypothesis for $\alpha$. Hence, $1 \leq k \leq m$ and $i\in \mathbb{N}$, $w\vdash_{\omega_1}\! \bigvee_{j\in \mathbb{N}} \neg\psi_{i, j}^{(k)}$. Therefore, $w\vdash_{\omega_1}\! \bigwedge_{k\leq m} \bigwedge_{i\in \mathbb{N}}  \bigvee_{j\in \mathbb{N}} \neg\psi_{i, j}^{(k)}$. Notice that $\vdash_{\omega_1}\! \bigwedge_{k\leq m} \bigwedge_{i\in \mathbb{N}}  \bigvee_{j\in \mathbb{N}} \neg\psi_{i, j}^{(k)} \leftrightarrow \neg \varphi$. Thus, $w \vdash_{\omega_1}\! \neg \varphi$. Moreover, the above proof shows that if  $w \vdash_{\omega_1}\! \neg \varphi$ then  $w \nvdash_{\omega_1}\!  \varphi$.		
		\item[2.] For the proof of left to right direction, suppose $w\cup\{\varphi\}$ is inconsistent in $\mathcal{H}_{\mathsf{DPL}_{\omega_1}}$. Then,  $w \vdash_{\omega_1}\! \neg \varphi$. But this implies that $w \nvdash_{\omega_1}\!  \varphi$, by Part 1. For the other direction,
		suppose $w\cup\{\varphi\}$ is consistent. Then, $w \nvdash_{\omega_1}\! \neg \varphi$. Hence, again by Part 1, $w \vdash_{\omega_1}\! \varphi$.
	\item[3.] This immediately follows from the fact that $w$ is saturated and therefore, is negation complete.
	\end{itemize}
	\vspace*{-4mm}
\end{proof}

\begin{coro}
Every saturated subset of $\mathbb{F}$ is consistent  in $\mathcal{H}_{\mathsf{DPL}_{\omega_1}}$.
\end{coro}
\subsection{The countable fragments of $\mathsf{DPL}_{\omega_1}$} \label{The countable fragments of IDPL}
In this subsection, we introduce the notion of a {\em countable fragment} of $\mathsf{DPL}_{\omega_1}$ and show that any countable fragment of  $\mathsf{DPL}_{\omega_1}$ has the strong completeness property with respect to the class $\mathcal{DMP}$ of all dynamic Markov processes.  The typical  proof of strong completeness theorem for a countable fragment $\mathbb{A}$ of $\mathsf{DPL}_{\omega_1}$ (e.g. \cite[Theorem 5.5]{Iko}) comprises of introducing a canonical model whose set of worlds consist of maximally consistent subsets of  $\mathbb{A}$. Hence, this canonical model depends on the fragment $\mathbb{A}$. However, we will prove that the canonical model introduced in Definition \ref{cpmodel1} of Subsection \ref{Strong completeness of DPL} for $\mathsf{DPL}$, serves not only as a canonical model for the finitary logic $\mathsf{DPL}$, but also  as a canonical model for each countable fragment $\mathbb{A}$ of  $\mathsf{DPL}_{\omega_1}$.

\begin{definition}
For every formula $\varphi$ of $\mathcal{L}_{\mathsf{DPL}_{\omega_1}}$, we  inductively define the {\em formal negation} of $\varphi$, denoted by $\sim\!\varphi$, as follows:
\begin{itemize}
\item[1.] $\sim \!p= \neg p$ for all $p \in \mathbb{P}$,
\item[2.] $\sim\!(\neg \varphi)= \varphi$,
\item[3.]  $\sim\! (\bigwedge_{i\in I} \varphi _i)= \bigvee_{i\in I} \sim\! \varphi _i $ for $|I|<\aleph_1$,
\item[4.]  $\sim\! (L_r \varphi)= \neg M_{1-r} \sim\! \varphi$, and
\item[5.]  $\sim\! (\bigcirc\varphi)=\bigcirc  \sim\! \varphi$,
\end{itemize}
\end{definition}

In the following, we bring the notion of a countable fragment of $\mathcal{L}_{\mathsf{DPL}_{\omega_1}}\!$, as it appears in \cite{Iko}. This notion stems in admissible model theory \cite{Bar, Keisler}.

\begin{definition}
We say that a (countable) set $\mathbb{A}$ of $\mathcal{L}_{\mathsf{DPL}_{\omega_1}}\!\!$-formulas is a {\em  (countable) fragment} if it satisfies the following closure properties:
\begin{itemize}
\item[1.] $\mathbb{F} \subseteq \mathbb{A}$;
\item[2.] $\mathbb{A}$ is closed under the finite Boolean connectives $\neg$ and $\wedge$,
\item[3.] $\mathbb{A}$ is closed under the probability and next operators, i.e., if $\varphi \in \mathbb{A}$ and $r\in \mathbb{Q} \cap [0,1]$, then $L_r \varphi \in \mathbb{A}$ and $\bigcirc\varphi \in \mathbb{A}$,
\item[4.] $\mathbb{A}$ is closed under the formal negation $\sim$, and
\item[5.] if $\varphi \in \mathbb{A}$ and $\psi$ is a subformula of $\varphi$, then $\psi \in \mathbb{A}$.
\end{itemize}
\end{definition}

We recall that the set $\text{Sub}(\varphi)$ of all subformulas of $\varphi$ can be defined as the smallest set of formulas of $\mathcal{L}_{\mathsf{DPL}_{\omega_1}}$ satisfying the following conditions:
\begin{itemize}
	\item[1.] $\varphi \in \textnormal{Sub}(\varphi)$;
	\item[2.] if $\star \psi \in \textnormal{Sub}(\varphi)$ for $\star\in\{\neg, L_r, \bigcirc\}$,  then $\psi \in \textnormal{Sub}(\varphi)$;
	\item[3.]  $\bigwedge_{i\in I} \psi_i \in \textnormal{Sub}(\varphi)$, then $\psi_i \in \textnormal{Sub}(\varphi)$ for each $i\in I$.
\end{itemize}

Notice that clearly $\mathbb{F}$ is a countable fragment of  $\mathcal{L}_{\mathsf{DPL}_{\omega_1}}$. Moreover, it is easy to see that for a given  (countable) set $\Gamma$ of formulas of $\mathcal{L}_{\mathsf{DPL}_{\omega_1}}\!\!$, there exists the smallest (countable) fragment $\mathbb{A}$ containing $\Gamma$. The following standard  lemma  which is usually called the Lindenbaum lemma, shows that for a given countable fragment $\mathbb{A}$, any consistent set $\Gamma$ of formulas in  $\mathbb{A}$ can be extended to a maximal consistent set $\Gamma^+\subseteq \mathbb{A}$. Here maximal is with respect to inclusion relation in set of all consistent subsets of $\mathbb{A}$ containing $\Gamma$.

\begin{lemma}[Lindenbaum lemma for countable fragments] \label{LinA}
 Let  $\mathbb{A}$ be a countable fragment of $\mathcal{L}_{\mathsf{DPL}_{\omega_1}}$ and
 $\Gamma\cup \{\varphi\} \subseteq\mathbb{A}$. If $\Gamma\nvdash_{\omega_1}\!\varphi$, then there exists a  maximal consistent subset $\Gamma^+$ of $\mathbb{A}$ such that $\Gamma \subseteq \Gamma^+$ and  $\Gamma^+  \nvdash_{\omega_1}\! \varphi$.
 \end{lemma}
\begin{proof}
Similar to  Lemma \ref{lin}, we should  inductively define a sequence
$\Gamma_0  \subseteq \Gamma_1  \subseteq  \dots \subseteq\Gamma_k  \subseteq \dots$ such that $\Gamma_k \nvdash_{\omega_1}\! \varphi$ for each $k \in \omega$.
The only new point here is how to use the countability of $\mathbb{A}$ to maintain $\Gamma_{k+1}$ from $\Gamma_k$.

First, suppose that $\varphi_0, \; \varphi_1,\;   \varphi_2, \;\dots$ and $s_0, \; s_1,\;  s_2, \;\dots$ are two enumerations of  all formulas of  $\mathbb{A}$ and all rational numbers in $[0, 1]$, respectively.  Then, take
 $\Gamma_0 := \Gamma$. Clearly,  $\Gamma_0 \nvdash_{\omega_1}\! \varphi$. Now, let $\Gamma_k$  be given such that $\Gamma_k \nvdash_{\omega_1}\! \varphi$. Then, by Part 3 of Lemma \ref{prop1}, either $\Gamma_k, \varphi_k \nvdash_{\omega_1}\! \varphi$ or $\Gamma_k, \neg\varphi_k \nvdash_{\omega_1}\! \varphi$. So we have two following cases 1 and 2:
\begin{itemize}
\item[1.] If $\Gamma_k, \varphi_k \nvdash_{\omega_1}\!\varphi$. Then, we define $\Gamma_{k+1} := \Gamma_k \cup \{\varphi_k\}$.
\item[2.] If  $\Gamma_k, \neg\varphi_k \nvdash_{\omega_1}\!\varphi$.  In this case, we define
	\begin{equation*}
	\Gamma_{k+1} := \left\{
	\begin{array}{rl}
	\Gamma_k \cup \{\neg\varphi_k, \neg \theta_i\}\;\;\;\;\;\;\;\;\;\;\;\;\;\;\;\;\;
	& \text{if } \varphi_k = \bigwedge_{i\in \mathbb{N}}\theta_i \text{ and $i\in \mathbb{N}$ is the least}\\
	& \text{number with }\Gamma_k, \neg\varphi_k, \neg\theta_i \nvdash_{\omega_1}\! \varphi,\\
\Gamma_k \cup \{\neg\varphi_k, \neg\bigcirc^n L_{s_l}\theta\}\;\;\;\;\;\;\;\!
	& \text{if }  \varphi_k = \bigcirc^n L_r\theta \;\text{and $l$ is the least number}\\
	& \text{such that $s_l<r$  and}\; \Gamma_k, \neg\varphi_k, \neg\bigcirc^n L_{s_l}\theta \nvdash_{\omega_1}\! \varphi,\\
\Gamma_k \cup \{\neg \varphi_k\}\;\;\;\;\;\;\;\;\;\;\;\;\;\;\;\;\;\;\;\;\;\;\;\;
	& \text{otherwise}\\
	\end{array} \right.
	\end{equation*}
\end{itemize}
where $r, s_l \in \mathbb{Q}\cap [0, 1]$ and $n\in \mathbb{N}$. Note that in Case 2 where $\varphi_k = \bigwedge_{i\in \mathbb{N}}\theta_i$, by Conj rule, there exists $i\in \mathbb{N}$ such that $\Gamma_k, \neg\varphi_k, \neg\theta_i \nvdash_{\omega_1}\! \varphi$ . Likewise, by Arch rule, there exists $l$ such that $s_l<r$  and $\Gamma_k, \neg\varphi_k, \neg\bigcirc^n L_{s_l}\theta \nvdash_{\omega_1}\! \varphi$, provided $\varphi_k = \bigcirc^n L_r\theta $. Hence, in either cases $\Gamma_{k+1}$ is well defined. Now, put
\begin{equation*}
\Gamma^+ := \bigcup_{k \in \omega}\Gamma_k.
\end{equation*}
Then, it can be shown that  $\Gamma^+$   satisfies all the desired properties.

First, note that all $\Gamma_k$'s are consistent in $\mathcal{H}_{\mathsf{DPL}_{\omega_1}}$, since $\Gamma_k \nvdash_{\omega_1}\! \varphi$. Next, by transfinite induction on the length of derivation of $\psi$ from $\Gamma^+$,  we can show that  $\Gamma^+$  is deductively closed, that is, for any formula $\psi \in \mathbb{A}$, if $\Gamma^+  \vdash_{\omega_1}\! \psi$ then $\psi \in \Gamma^+$. So $\Gamma^+$ is consistent.
 Furthermore, by the above construction, $\Gamma^+$ is negation complete. Thus, $\Gamma^+$ is maximal. Finally, since $\Gamma^+$ is deductively closed, it follows that $\Gamma^+  \nvdash_{\omega_1}\! \varphi$.
\end{proof}

\begin{coro} \label{LinB}
	Let  $\mathbb{A}$ be a countable fragment of $\mathcal{L}_{\mathsf{DPL}_{\omega_1}}$ and
	$\Gamma\cup \{\varphi\} \subseteq\mathbb{A}$. If $\Gamma\nvdash_{\omega_1}\!\varphi$, then there exists a   saturated subset $w$ of $\mathbb{F}$ such that $\Gamma\cup w \nvdash_{\omega_1}\!\varphi$.
\end{coro}
\begin{proof}
Let $\Gamma^+$ be the maximal consistent subset of $\mathbb{A}$ containing
$\Gamma$ and $\Gamma^+ \nvdash_{\omega_1}\!\varphi$, obtained from Lemma \ref{LinA}. Consider $w := \Gamma^+ \cap \mathbb{F}$. It can be easily check that $w$ is a saturated subset of $\mathbb{F}$ and $\Gamma\cup w  \nvdash_{\omega_1} \! \varphi$.
\end{proof}

The following lemma gives a criterion for satisfaction of a formula $\varphi \in \mathcal{L}_{\mathsf{DPL}_{\omega_1}}$ in the canonical model $\mathfrak{M}_c$, introduced in Definition \ref{cpmodel1}.

\begin{lemma}[Truth lemma] \label{truth}
Consider the canonical model $\mathfrak{M}_c$. Then, for any formula $\varphi \in \mathcal{L}_{\mathsf{DPL}_{\omega_1}}$ and $w \in \Omega_c$, we have
\begin{itemize}
\item[1.]  $\mathfrak{M}_c, w\vDash \varphi \text{ iff } w\vdash_{\omega_1}\! \varphi$;
\item[2.] the set $[\![\varphi]\!]_{\mathfrak{M}_c} = \{w\in \Omega_c\;|\;\mathfrak{M}_c, w \vDash \varphi\} $  is measurable.
\end{itemize}
\end{lemma}
\begin{proof}
\begin{itemize}
	\item[1.] We first show the statement for $\varphi\in \Delta_\infty$. We use the transfinite induction to verify that  1 holds for each formula $\varphi \in \Delta_\alpha$ and  $\alpha < \nu$. Notice that the claim is true for  $\varphi \in \Delta_0$, since
	in this case,  $w \vdash_{\omega_1}\! \varphi$ implies $\varphi\in w$. So by  Lemma \ref{truth2},  the induction base is established.
	
	Now, let $\varphi \in \Delta_{\alpha+1}$. Then by Part 4 of Remark \ref{repre1},  $\varphi$ is provably equivalent with a formula of the form $\bigvee_{k=1}^{m} \gamma_k$ such that  for each  $1 \leq k \leq m$, $\gamma_k= \bigvee_{i\in \mathbb{N}} \bigwedge_{j\in \mathbb{N}} \psi_{i,j}^{(k)}$ with $\psi_{i, j}^{(k)}\in \Delta_\alpha$.
	If $\mathfrak{M}_c, w\vDash \varphi$, then there are $1 \leq k \leq m$ and $i \in \mathbb{N}$ such that for each $j\in \mathbb{N}$, we have $\mathfrak{M}_c, w\vDash \psi_{i,j}^{(k)}$. So by the induction hypothesis, there are $1 \leq k \leq m$ and $i \in \mathbb{N}$ such that $w\vdash_{\omega_1}\!  \psi_{i,j}^{(k)}$, for each $j\in \mathbb{N}$. Hence, $w\vdash_{\omega_1}\! \varphi$. The other direction can be proved similarly.
	
	Now, if $\varphi$ is an arbitrary formula in  $\mathcal{L}_{\mathsf{DPL}_{\omega_1}}$, then by Lemma \ref{equi}, there exists $\varphi'\in  \Delta_\infty$ such that
	$\vdash_{\omega_1}\! \varphi\leftrightarrow\varphi'$. So by Theorem \ref{sound2}, we have $\mathfrak{M}_c, w\vDash\varphi\leftrightarrow\varphi'$ for each $w \in \Omega_c$.
	This implies that
	$\mathfrak{M}_c, w\vDash\varphi$ iff  $\mathfrak{M}_c, w\vDash \varphi'$ for each $w \in \Omega_c$. Furthermore, $w\vdash_{\omega_1}\! \varphi$ iff $w\vdash_{\omega_1}\! \varphi'$ for each $w \in \Omega_c$. Hence, the claim is easy to verify.
	\item[2.] Similar to 1, we use the transfinite induction and assume that $\varphi\in \Delta_\infty$. The measurability of $[\![\varphi]\!]_{\mathfrak{M}_c}$
	for $\varphi\in \Delta_0$ can be concluded from Lemma \ref{truth2}. Now assume that $\varphi \in \Delta_{\alpha+1}$ is of the form $\bigvee_{k\leq m} \gamma_k$ such that  for each  $1 \leq k \leq m$, $\gamma_k= \bigvee_{i\in \mathbb{N}} \bigwedge_{j\in \mathbb{N}} \psi_{i,j}^{(k)}$ with $\psi_{i, j}^{(k)}\in \Delta_\alpha$. 	Then, the induction hypothesis holds for $\varphi$, since
	\[ [\![\varphi]\!]_{\mathfrak{M}_c} = \bigcup_{k\leq m}  \;\bigcup_{i\in \mathbb{N}} \;\bigcap_{j\in \mathbb{N}} \; [\![\psi_{i, j}^{(k)} ]\!]_{\mathfrak{M}_c}. \] 	
\end{itemize}
	\vspace*{-4mm}
\end{proof}

\begin{theo} [Strong completeness of countable fragments] \label{com2}
	Let $\mathbb{A}$ be a countable fragment of $\mathcal{L}_{\mathsf{DPL}_{\omega_1}}$ and $\Gamma\cup \{\varphi\}$ be a set of formulas of $\mathbb{A}$. Then, $\Gamma \vDash \varphi$ implies $\Gamma \vdash_{\omega_1}\! \varphi$.
\end{theo}
\begin{proof}
	Suppose that $\Gamma  \nvdash_{\omega_1}\!\varphi$. Then, by Lemma \ref{LinA} and Corollary \ref{LinB}, there exist a maximal consistent $\Gamma^+$ of $\mathbb{A}$
	and a saturated set $w\in \Omega_c$ such that $\Gamma \cup w \subseteq \Gamma^+$
	and $\Gamma^+ \nvdash_{\omega_1}\! \varphi$. Notice that for each $\theta \in \Gamma^+$, $w\cup \{\theta\}$ is consistent. Hence, by Part 2 of Lemma \ref{negcompleteF},  $w  \vdash_{\omega_1}\!\theta$, for each $\theta \in \Gamma^+$. On the other hand, $w  \vdash_{\omega_1}\!\varphi $. So in the light of
 Part 1 of Lemma \ref{truth}, $\mathfrak{M}_c, w \vDash \Gamma$, while $\mathfrak{M}_c, w \nvDash \varphi$. Thus, we have $\Gamma \nvDash \varphi$.
\end{proof}

\begin{coro} [Weak completeness of $\mathcal{H}_{\mathsf{DPL}_{\omega_1}}$] \label{com1}
	Let  $\varphi \in \mathcal{L}_{\mathsf{DPL}_{\omega_1}}$ with $\nvdash_{\omega_1}\!\varphi$, Then there exists $w \in \Omega_c$ such that  $\mathfrak{M}_c, w\vDash \varphi$. Hence, $\varphi$ is not a tautology.
\end{coro}
\begin{proof}
	Let $\nvdash_{\omega_1}\!\varphi$. Then, find a countable fragment of  $\mathcal{L}_{\mathsf{DPL}_{\omega_1}}$ containing $\varphi$. Then, we can apply Theorem \ref{com2}, and find a saturated set $w$ such that $w\nvdash_{\omega_1}\!\varphi$ and $\mathfrak{M}_c, w\nvDash \varphi$. Thus, $\varphi$ is not a tautology.
\end{proof}

\section{Definability} \label{Definability}
In this section, we turn our attention to studying (frame) definability of certain dynamic and probabilistic properties of dynamic Markov processes.  We will show that some dynamic properties such as measure-preserving, ergodicity and mixing are definable within $\mathsf{DPL}$ and $\mathsf{DPL}_{\omega_1}$.  Moreover, we consider the infinitary probability logic with initial distribution ($\mathsf{InPL}_{\omega_1}$) by disregarding the dynamic operator. We illustrate that the strong expressive power of ${\mathsf{InPL}_{\omega_1}}$ would allow us to define $n$-step transition probabilities $T^n$ of Markov kernel $T$. This implies that some probabilistic properties of Markov processes such as stationary, invariance, irreducibility and recurrence are definable within ${\mathsf{InPL}_{\omega_1}}$.

\begin{definition}\label{def}
Let $\mathcal{C} \subseteq \mathcal{D}$ be two classes of dynamic Markov processes.
We say that the class $\mathcal{C}$ is {\em  $\mathsf{DPL}$-definable} (resp. {\em $\mathsf{DPL}_{\omega_1}$-definable}) {\em within $\mathcal{D}$}
if there exists a set  $\Gamma$  of formulas of $\mathcal{L}_{\mathsf{DPL}}$ (resp. $ \mathcal{L}_{\mathsf{DPL}_{\omega_1}}$) such that for each $\mathfrak{P} \in \mathcal{D}$, we have  $\mathfrak{P} \in \mathcal{C}$ if and only if $\mathfrak{P}\vDash \Gamma$.  A class  $\mathcal{C}$ is said to be {\em $\mathsf{DPL}$-definable} (resp. {\em $\mathsf{DPL}_{\omega_1}$-definable}) if it is  $\mathsf{DPL}$-definable (resp. $\mathsf{DPL}_{\omega_1}$-definable)  within the class $\mathcal{DMP}$ of all dynamic Markov processes.
\end{definition}

\paragraph{Measure-preserving functions.}
Recall that a measurable function $f$ on a probability space $\langle{\Omega, \mathcal{A}, \mu}\rangle$ is {\em measure-preserving} if $\mu( f^{-1}(A)) = \mu(A)$ for any $A\in \mathcal{A}$; see, e.g. \cite[Definition 2.1.1]{Douc}. In this case, we call the tuple $\langle{\Omega, \mathcal{A}, \mu, f}\rangle$ an {\em abstract dynamical system}.

\begin{pro} \label{Definability1}
The class $\mathcal{ADS}$ of all abstract dynamical system is $\mathsf{DPL}$-definable within the class $\mathcal{DPS}$ of all dynamic probability spaces by the set  \[\mathbf{M}:= \{L_r \bigcirc p\leftrightarrow\bigcirc L_r p\;|\; r\in \mathbb{Q}\cap[0, 1]\}.\]
\end{pro}
\begin{proof}
	Let $\mathfrak{P}= \langle{\Omega, \mathcal{A}, T, f}\rangle$ be a dynamic probability space. Then, by Part 4 of Definition \ref{frames}, there exists a probability measure $\mu:\mathcal{A}\to[0, 1]$ such that  $T(w, A)=\mu(A)$ for each $w\in \Omega$ and $A\in \mathcal{A}$. We must show that
	\[\mathfrak{P}\vDash \mathbf{M} \text{  if and only if }  f \text{ is measure-preserving}.\]
	First suppose that $f$ is not measure-preserving. Then, there exists $A\in \mathcal{A}$ such that $\mu( f^{-1}(A)) \neq \mu(A)$, say  $\mu( f^{-1}(A)) > \mu(A)$. Therefore,  there exists $r\in \mathbb{Q} \cap [0,1]$ such that	$T(w,  f^{-1}(A))\geq r$, but $ T(f(w), A) < r$ for some (equivalently, each) $w\in \Omega$. Now consider the valuation $v$ on  $\mathfrak{P}$ such that $v(p) = A$. In this case, we have $\mathfrak{M}, w\vDash  L_r  \bigcirc p$, but $\mathfrak{M}, w\nvDash \bigcirc  L_r  p$.  This means that  $\mathfrak{M}, w\nvDash L_r \bigcirc  p \to \bigcirc  L_r p$, a contradiction. The case $\mu( f^{-1}(A)) < \mu(A)$ can be proved similarly.
	\\
	Conversely,  suppose that $f$ is measure-preserving and $\mathfrak{M}$ is a model based on $\mathfrak{P}$ and $w\in \Omega$ such that $\mathfrak{M}, w\vDash L_r \bigcirc p$. Then, we have $T(w, [\![\bigcirc p]\!]_{\mathfrak{M}}) \geq r$. So $T(w, f^{-1}([\![p]\!]_{\mathfrak{M}})) \geq r$. Now since $f$ is measure-preserving, it follows that $T(f(w), [\![p]\!]_{\mathfrak{M}}) \geq r$. Thus, $\mathfrak{M}, f(w)\vDash L_r p$, and hence $\mathfrak{M}, w\vDash \bigcirc L_r p$. So, we have $\mathfrak{M}, w\vDash L_r \bigcirc p \to  \bigcirc  L_r p$.  Similarly, it can be shown that $\mathfrak{M}, w\vDash \bigcirc L_r p\to L_r \bigcirc p$.
\end{proof}

\paragraph{Ergodicity.}
Let $\langle{\Omega, \mathcal{A}, \mu, f}\rangle$  be a dynamic probability space. A set $A\in \mathcal{A}$ is {\em invariant} whenever $f^{-1}(A) = A$. The function $f$ is said to be {\em ergodic} if it has trivial invariant sets, that is, if $A\in \mathcal{A}$ is invariant, then either $\mu(A)= 0$ or $\mu(A)= 1$; see, e.g. \cite[Definition 5.1.10]{Douc}. In this case, the tuple $\langle{\Omega, \mathcal{A}, \mu, f}\rangle$ is called an {\em ergodic dynamical system}.

\begin{pro}
The class of all ergodic dynamical systems is $\mathsf{DPL}$-definable within the class of all dynamic probability spaces by the formula $\mathbf{E}:= (\bigcirc p \leftrightarrow p) \to (L_0 p \vee  L_1 p)$.
\end{pro}
\begin{proof}
	Suppose that $\mathfrak{P}= \langle{\Omega, \mathcal{A}, T, f}\rangle$ is a dynamic probability space. Then, there exists a probability measure $\mu:\mathcal{A}\to[0, 1]$ such that  $T(w, A)=\mu(A)$ for each $w\in \Omega$ and $A\in \mathcal{A}$. We need to show that
\[\mathfrak{P}\vDash \mathbf{E} \text{ if and only if } f \text{ is ergodic}.\]
	We first prove the left to right direction. Suppose $f$ is not ergodic. Then, there exist  $w\in \Omega$ and  $A\in \mathcal{A}$ such that  $f^{-1}(A) = A$ and neither $T(w, A)= 0$ nor $T(w, A)= 1$.
	Now consider the valuation $v$ on  $\mathfrak{P}$ such that $v(p) = A$. In this case,  $\mathfrak{M}, w\vDash  \bigcirc p \leftrightarrow p$, but $\mathfrak{M}, w\nvDash L_0 p \vee  L_1 p$. So  $\mathfrak{M}, w\nvDash \mathbf{E}$. The other direction can be shown in a similar way.
\end{proof}

\begin{coro}
The class of all ergodic abstract dynamical systems is $\mathsf{DPL}$-definable within the class of all dynamic probability spaces.
\end{coro}

\paragraph{Mixing.}
Let $\langle{\Omega, \mathcal{A}, \mu, f}\rangle$ be an abstract dynamical system. Then, $f$ is said to be ({\em strongly}) {\em mixing} if for all $A, B\in \mathcal{A}$, we have
\[\lim_{k\to\infty} \mu(f^{-k}(A)\cap B) = \mu(A)\mu(B).\]
See, e.g. \cite[Definition 3.4.1]{Dajani}. In this case, the tuple $\langle{\Omega, \mathcal{A}, \mu, f}\rangle$ is called a  ({\em strongly}) {\em mixing dynamical system}.

\begin{pro}
	The class of all mixing dynamical systems is $\mathsf{DPL}_{\omega_1}$-definable within the class of all abstract dynamical systems by the formula
	\[\mathbf{Mix}:= \mathbf{G}\; \wedge \;\mathbf{L}\]
	 where
\[\mathbf{G}:=  \bigwedge_{r, s\in \mathbb{Q} \cap [0, 1]} \big((L_{r} p \wedge  L_{s} q) \to \bigwedge_{0	<n\in \mathbb{N}} \; \bigvee_{m\in \mathbb{N}} \; \bigwedge_{k\geq m}\; L_{rs\dot{-}\frac{1}{n}}(\bigcirc^k p \wedge q) \big),\]
and
\[\mathbf{L}:=  \bigwedge_{r, s\in \mathbb{Q} \cap [0, 1]} \big((
M_{r} p \wedge  M_{s} q) \to \bigwedge_{0<n\in \mathbb{N}} \; \bigvee_{m\in \mathbb{N}} \;  \bigwedge_{k\geq m}\; M_{rs +\frac{1}{n}}(\bigcirc^k p \wedge q) \big).\]
\end{pro}
\begin{proof}
	Assume that $\mathfrak{P}= \langle{\Omega, \mathcal{A},\mu, f}\rangle$ is an abstract dynamical system. It is not hard to check that $\mathfrak{P}\vDash \mathbf{G}$ if and only if
	\[\lim_{k\to\infty} \mu(f^{-k}(A)\cap B) \geq \mu(A)\mu(B),\]
	for each  $A, B \in \mathcal{A}$. Likewise,
	$\mathfrak{P}\vDash \mathbf{L}$ if and only if
	\[\lim_{k\to\infty} \mu(f^{-k}(A)\cap B) \leq \mu(A)\mu(B),\]
	for each  $A, B \in \mathcal{A}$. Hence, we have $\mathfrak{P}\vDash \mathbf{Mix} \text{ if and only if } f \text{ is mixing}$.
\end{proof}
\subsection{Definability of $n$-step transition probabilities} \label{nstep}
Here, we consider pure probabilistic dynamic Markov processes, or simply,  by Remark \ref{P}, Markov processes. That means the dynamic part is disregarded. Hence, we denote  the resulting finitary and infinitary logics by $\mathsf{PL}$ and $\mathsf{PL}_{\omega_1}$, respectively. This assumption  is crucially for convenience, since the dynamic part  does not play a key role in this situation. Further, this part seems a new contribution for probability logic.

In order to make our studies as pervasive as possible, we focus on Markov processes with initial distributions. So,  we extend the ordinary language of probability logic $\mathcal{L}_{\mathsf{PL}}$ by some probability operators  $\mathsf{L}_r^{\langle{0}\rangle}$ for each $r\in \mathbb{Q} \cap [0,1]$.

\begin{definition}
 Let $\mathcal{L}_{\mathsf{PL}}$ be the language of the usual probability logic containing the probability operators $L_r$ for each $r\in \mathbb{Q} \cap [0,1]$. We enrich this language by adding some probability operators $\mathsf{L}_r^{\langle{0}\rangle}$ for each $r\in \mathbb{Q} \cap [0,1]$. These operators  are called {\em initial probability operators}.  We denote the expanding language by
 $\mathcal{L}_{\mathsf{InPL}}$. If  $\mathbb{P}$  is a countable set of propositional variables, then the formulas of $\mathcal{L}_{\mathsf{InPL}}$ can be recursively defined as follows:
\begin{align*}
	\varphi ::&=\ \ p\mid \neg\varphi \mid \varphi \land \varphi\mid L_r \varphi  \mid \mathsf{L}_r^{\langle{0}\rangle} \varphi
\end{align*}
where $p \in \mathbb{P}$ and $r\in \mathbb{Q} \cap [0,1]$. In the sequel, the language $\mathcal{L}_{\mathsf{InPL}_{\omega_1}}$  is the infinitary logic obtained from $\mathcal{L}_{\mathsf{InPL}}$ in the similar way as $\mathcal{L}_{\mathsf{DPL}_{\omega_1}}$.
\end{definition}

\begin{definition}
A quadruple $\langle{\Omega, \mathcal{A}, \pi, T}\rangle$ is called a {\em Markov process with initial distribution}  whenever $\langle \Omega, \mathcal{A}, T\rangle$ is a Markov process and $\pi:\mathcal{A} \to [0, 1]$ is a $\sigma$-additive probability measure. If $v: \mathbb{P} \to \mathcal{A}$ is a valuation function, then the resulting tuple $ \langle{\Omega, \mathcal{A}, \pi, T, v}\rangle$ is a {\em Markov model with initial distribution}, abbreviated as  {\em In-model}.
\end{definition}

Given an In-model $\mathfrak{M} = \langle{\Omega, \mathcal{A}, \pi, T, v}\rangle$, the satisfiability relation for the new probability formulas of $ \mathsf{L}_r^{\langle{0}\rangle} \varphi$ in a world $w\in\Omega$ is inductively defined as:
\begin{itemize}
	\item[] $\mathfrak{M}, w \vDash \mathsf{L}_r^{\langle{0}\rangle} \varphi \;$ iff $\; \pi([\![\varphi]\!]_{\mathfrak{M}}) \geq r$.
\end{itemize}

One can consider the logic $\mathsf{InPL}$ (resp.  ${\mathsf{InPL}_{\omega_1}}$) to be the set of all valid formulas of $\mathcal{L}_{\mathsf{InPL}}$ (resp.  $\mathcal{L}_{\mathsf{InPL}_{\omega_1}}$) over the class of all Markov processes with initial distribution.  Subsequently, one can regard the proof system $\mathcal{H}_{\mathsf{InPL}}$ (resp.  $\mathcal{H}_{\mathsf{InPL}_{\omega_1}}$) containing the usual axioms and rules of probability logic for the expanding language $\mathcal{L}_{\mathsf{InPL}}$ (resp.  $\mathcal{L}_{\mathsf{InPL}_{\omega_1}}$)  together with Harsanyi axiom schemes of the form $H_1$ and $H_2$ for probability operators $\mathsf{L}_r^{\langle{0}\rangle}$, introduced in Definition \ref{Axiom}, that is
$\mathsf{L}_r^{\langle{0}\rangle} \varphi \to \mathsf{L}_1^{\langle{0}\rangle}\mathsf{L}_r^{\langle{0}\rangle} \varphi$
and $\neg \mathsf{L}_r^{\langle{0}\rangle} \varphi \to \mathsf{L}_1^{\langle{0}\rangle} \neg\mathsf{L}_r^{\langle{0}\rangle} \varphi$, respectively.
Once these proof systems are set, one can easily extend Part 2 of Theorem \ref{MP and Pure} and Theorem \ref{PS}, showing the following results.

\begin{theo}
\begin{itemize}
\item[1.] The proof system $\mathcal{H}_{\mathsf{InPL}}$ is strongly sound and complete for the class $\mathcal{I}n\mathcal{MP}$ of all Markov processes with initial distribution.
\item[2.] The proof system $\mathcal{H}_{\mathsf{InPL}_{\omega_1}}$  is weakly sound and complete for the class $\mathcal{I}n\mathcal{MP}$.
\end{itemize}
\end{theo}

Notice that the both notions of  {\em $\mathsf{InPL}$-definability} and {\em $\mathsf{InPL}_{\omega_1}$-definability} can be defined the same way as Definition \ref{def}. Now having said that, we verify that the expressive power of $\mathcal{L}_{\mathsf{InPL}_{\omega_1}}$ is strong enough for showing that the {\em $n$-step transition probability $T^n$} of Markov kernel $T$, is $\mathsf{InPL}_{\omega_1}$-definable. This fact subsequently implies that many natural stochastic properties of Markov processes such as stationary, invariance, irreducibility and recurrence are $\mathsf{InPL}_{\omega_1}$-definable. In the following, let us first recall the notion of $n$-step transition probability $T^n$ of Markov kernel $T$, for each $n\in \mathbb{N}$.

\begin{definition} \cite[Subsection 1.2.1]{Meyn} \label{ntimes}
	Let $\langle{\Omega, \mathcal{A}, \pi, T}\rangle$ be a Markov process with initial distribution. 
	Then, for each $n \in \mathbb{N}$, the {\em $n$-step transition probability $T^n$} is defined as follows. Let $T^0:= \pi$, $T^1:= T$ and $T^{n+1}$ for each $n\geq 1$,
	\[T^{n+1}(w, A):=\int_{\Omega} T^n (y, A)\;T(w, dy)\] 
	for every $w\in \Omega$ and $A\in \mathcal{A}$. Notice that here it is assumed for every $A\in \mathcal{A}$, $T^{n} (., A)$ is a measurable function on $\mathcal{A}$ and that the right hand of the above equation is just integration of the function $T^{n} (., A)$ with respect to measure $T(w, .)$.
\end{definition}

It is known that for each $n\in \mathbb{N}$, $T^n: \Omega\times \mathcal{A} \to [0, 1]$ are Markov kernels, see Definition \ref{Markov}. Towards defining  $n$-step transition probability $T^n$, we introduce the infinitary formulas $\mathbb{L}_r^{\langle{n}\rangle} (\varphi)$ and $\mathbb{M}_r^{\langle{n}\rangle} (\varphi)$ for each $\varphi \in \mathcal{L}_{\mathsf{InPL}_{\omega_1}}$,  $n\in \mathbb{N}$ and $r\in \mathbb{Q} \cap [0,1]$.

\begin{definition}
For given $\varphi \in \mathcal{L}_{\mathsf{InPL}_{\omega_1}}$, $n\in \mathbb{N}$ and $r\in \mathbb{Q} \cap [0,1]$, we define the infinitary formulas $\mathbb{L}_r^{\langle{n}\rangle} (\varphi)$ and $\mathbb{M}_r^{\langle{n}\rangle} (\varphi)$ inductively on $n$ as follows:
\begin{itemize}
\item Let $\mathbb{L}_r^{\langle{0}\rangle} (\varphi) := \mathsf{L}_r^{\langle{0}\rangle}  \varphi$, $\mathbb{M}_r^{\langle{0}\rangle} (\varphi) := \mathsf{M}_r^{\langle{0}\rangle}  \varphi = \mathsf{L}_{1-r}^{\langle{0}\rangle}  \neg\varphi$, and that  $\mathbb{L}_r^{\langle{1}\rangle} (\varphi) := L_r \varphi$ and $\mathbb{M}_r^{\langle{1}\rangle} (\varphi) := M_r \varphi= L_{1-r} \neg \varphi$, for each $r\in \mathbb{Q} \cap [0,1]$;
\item Suppose that for each $\psi \in \mathcal{L}_{\mathsf{InPL}_{\omega_1}}$ and $s\in \mathbb{Q} \cap [0,1]$, the formulas $\mathbb{L}_s^{\langle{n}\rangle} (\psi)$ and  $\mathbb{M}_s^{\langle{n}\rangle} (\psi)$ are already defined in $\mathcal{L}_{\mathsf{InPL}_{\omega_1}}$. Then, set
\[\mathbb{L}_r^{\langle{n+1}\rangle} (\varphi) :=   \bigwedge_{l\geq 1} \;\;\bigvee_{k\geq 1} \;\;\bigvee_{r_0+\dots+r_{k-1}\geq r\dot{-}\frac{1}{l}} \;\;\bigwedge_{i= 0}^{k-1}\; \mathbb{L}_{\frac{kr_i}{i}}^{\langle{n}\rangle} \big(\mathbb{L}_{\frac{i}{k}}^{\langle{n}\rangle}  (\varphi)  \wedge \neg \mathbb{L}_{\frac{i+1}{k}}^{\langle{n}\rangle} (\varphi)  \big),\]
where $r_0, \dots, r_{k-1} \in \mathbb{Q} \cap [0,1]$. Respectively, consider
$\mathbb{M}_r^{\langle{n+1}\rangle} (\varphi) := \mathbb{L}_{1-r}^{\langle{n+1}\rangle} (\neg \varphi)$.
\end{itemize}
\end{definition}

The following lemma exhibits the connection of the infinitary formulas
$\mathbb{L}_r^{\langle{n}\rangle} (\varphi)$ and the transition probabilities $T^n$, for all  $\varphi\in \mathcal{L}_{\mathsf{InPL}_{\omega_1}}$, $n\in \mathbb{N}$ and $r\in \mathbb{Q} \cap [0,1]$.

\begin{lemma} \label{semequ}
Suppose that  $\mathfrak{M} = \langle{\Omega, \mathcal{A}, \pi, T,  v}\rangle$ is an In-model and $w \in \Omega$. Then, for each $\varphi\in \mathcal{L}_{\mathsf{InPL}_{\omega_1}}$, $n\in \mathbb{N}$ and $r\in \mathbb{Q} \cap [0,1]$, we have
\[\mathfrak{M}, w \vDash \mathbb{L}_r^{\langle{n}\rangle} (\varphi) \; \text{ iff }\; T^n(w, [\![\varphi]\!]_{\mathfrak{M}}) \geq r.\]
\end{lemma}
\begin{proof}
This is proved by induction on $n$. The cases of $n=0 $ and $n=1$ hold by definition. Assume that the equivalence holds for $n\geq 1$,  for each $s\in \mathbb{Q} \cap [0,1]$ and $\psi\in \mathcal{L}_{\mathsf{InPL}_{\omega_1}}$. Now to show the equivalence for $n+1$, let  $\varphi\in \mathcal{L}_{\mathsf{InPL}_{\omega_1}}$ and $r\in \mathbb{Q} \cap [0,1]$. Put  $A =[\![\varphi]\!]_{\mathfrak{M}}$ for some $A\in \mathcal{A}$ and assume that
\[T^{n+1}(w, A):=\int_{\Omega}  T^n(y, A)\; T(w, dy)\geq r.\]
For $k\geq 1$, take the partition $\{0, \dfrac{1}{k}, \dots, \dfrac{k-1}{k}, 1\}$ from $[0, 1]$. Then, by definition of Lebesgue integration (e.g., see \cite{Dudley}), we have
\[\int_{\Omega} T^n(y, A) \; T(w, dy)= \sup\big\{\sum_{i=0}^{k-1}\; \dfrac{i}{k} \;T^n(w, A_i)\;|\; k\geq 1 \big\} \geq r\]
where for each $0\leq i\leq k-1$, we have
$A_i= \{y\in \Omega\;|\; \dfrac{i}{k} \leq T^n(y, A) < \dfrac{i+1}{k} \}.$
Note that by induction hypotheses, for each $0\leq i\leq k-1$, we obtain that
\[ y\in A_i \;\;\text{ iff } \;\;\mathfrak{M}, y \vDash \mathbb{L}_{\frac{i}{k}}^{\langle{n}\rangle}  (\varphi)  \wedge \neg\mathbb{L}_{\frac{i+1}{k}}^{\langle{n}\rangle}  (\varphi) ,\]
 for each $y\in \Omega$. This means  $A_i= [\![\mathbb{L}_{\frac{i}{k}}^{\langle{n}\rangle}  (\varphi) \wedge \neg \mathbb{L}_{\frac{i+1}{k}}^{\langle{n}\rangle}  (\varphi) ]\!]_{\mathfrak{M}} $.  Let $\alpha_k = \sum_{i=0}^{k-1}\; \dfrac{i}{k} \;T^n(w, A_i)$. Then, $\sup \;\!\{\alpha_k\;| \; k\geq1\} \geq r$  is equivalent to this fact that for each natural number $l\geq 1$ there exists $k \geq 1$ such that $\alpha_k \geq r \dot{-} \frac{1}{l}$. But,  $\alpha_k \geq r \dot{-} \frac{1}{l}$ is the same as existence of the rational numbers $r_0, \dots, r_{k-1} \in \mathbb{Q} \cap [0,1]$ such that $r_0+\dots+r_{k-1}\geq  r \dot{-} \frac{1}{l}$ and $T^n(w, A_i) \geq\dfrac{k r_i}{i}$
 for each  $0\leq i\leq k-1$. Notice that, by induction hypothesis, the latter inequality is equivalent to
  \[\mathfrak{M}, w \vDash  \mathbb{L}_{\frac{kr_i}{i}}^{\langle{n}\rangle} \big(\mathbb{L}_{\frac{i}{k}}^{\langle{n}\rangle}  (\varphi)  \wedge \neg \mathbb{L}_{\frac{i+1}{k}}^{\langle{n}\rangle} (\varphi)  \big).\]
 Thus,  we have
 \[T^{n+1}(w, [\![\varphi]\!]_{\mathfrak{M}}) \geq r  \;\;\;\text{iff}\;\;\;  \mathfrak{M}, w \vDash   \bigwedge_{l\geq 1} \;\;\bigvee_{k\geq 1} \;\;\bigvee_{r_0+\dots+r_{k-1}\geq r \dot{-} \frac{1}{l}} \;\;\bigwedge_{i= 0}^{k-1}\; \mathbb{L}_{\frac{kr_i}{i}}^{\langle{n}\rangle} \big(\mathbb{L}_{\frac{i}{k}}^{\langle{n}\rangle}  (\varphi)  \wedge \neg \mathbb{L}_{\frac{i+1}{k}}^{\langle{n}\rangle} (\varphi)  \big).\]
 Hence,
  \[T^{n+1}(w, [\![\varphi]\!]_{\mathfrak{M}}) \geq r  \;\;\;\text{iff}\;\;\;  \mathfrak{M}, w \vDash   \mathbb{L}_{r}^{\langle{n+1}\rangle} (\varphi).\]
\end{proof}

The remainder of this section is devoted to exploring the $\mathsf{InPL}_{\omega_1}$-definability of certain stochastic properties of Markov processes with initial distribution. For more information about such properties, we refer the reader to any textbook on Markov processes; see, e.g., \cite{Douc, Meyn}.

\paragraph{Stationary.}
Let $\langle{\Omega, \mathcal{A}, T}\rangle$ be a Markov process. A probability measure $\mu$ on $\mathcal{A}$ is said to be {\em invariant} with respect to $T$ if $\mu T= \mu$, i.e.,
\[\mu T (A) = \int_\Omega T(w, A)\; \mu(dw)\]
for all $A\in \mathcal{A}$. A  Markov process  with initial distribution $\langle{\Omega, \mathcal{A}, \pi, T}\rangle$ is {\em stationary} whenever $\pi$ is invariant with respect to $T$.  

\begin{pro}
The class of all stationary Markov processes with initial distribution is $\mathsf{InPL}_{\omega_1}$-definable by the formula
	\[\mathbf{S}:=  \bigwedge_{r\in \mathbb{Q} \cap [0,1]}\Big(\mathsf{L}_r^{\langle{0}\rangle} (p) \leftrightarrow  \bigvee_{l\geq 1} \;\;\bigwedge_{k  \geq 1}\; \bigvee_{r_0+\dots+r_{k-1}\geq r\dot{-}\frac{1}{l}}\; \bigwedge_{i= 0}^{k-1} \mathsf{L}^{\langle{0}\rangle}_{\frac{kr_i}{i}}\big(L_{\frac{i}{k}} p \wedge \neg L_{\frac{i+1}{k}} p \big)\Big).\]
\end{pro}
\begin{proof}
	The proof is similar to the proof of Lemma \ref{semequ}. Let $\mathfrak{P}= \langle{\Omega, \mathcal{A}, \pi, T}\rangle$ be a Markov process with initial distribution. First notice that  we have
	\[\int_{\Omega}  T(w, A)\; \pi(dw)= \sup\big\{\sum_{i=0}^{k-1}\; \dfrac{i}{k} \;\pi(A_i)\;|\; k\geq 1\big\}\]
	where for each $0\leq i\leq k-1$,  $A_i= \{y\in \Omega\;|\; \dfrac{i}{k} \leq T(y, A) < \dfrac{i+1}{k} \}$.
	Moreover,  for every  $A \in \mathcal{A}$, $\pi T(A)= \pi(A)$ if and only if  we have
\begin{align*}
\pi (A) \geq r \;\; \text{ iff }\;\;  \pi T(A) \geq r \;\; & \text{ iff }\;\;\sup\big\{\sum_{i=0}^{k-1}\; \dfrac{i}{k} \;\pi(A_i)\;|\; k\geq 1 \big\} \geq r,
\end{align*}
for each $r \in \mathbb{Q} \cap [0,1]$.	 So it is not hard to see that $\mathfrak{P}\vDash \mathbf{S}\; \text{ iff }\;  \mathfrak{P} \text{ is stationary}.$
\end{proof}

\paragraph{Irreducibility.} Let $\langle{\Omega, \mathcal{A}, \pi, T}\rangle$ be a Markov process with initial distribution. We say that this process is {\em irreducible} if for all $A\in \mathcal{A}$ and  $w\in \Omega$, whenever $\pi (A)>0$, there exists  $n \geq 1$ such that  $T^n(w, A)>0$, see  \cite[Definition 9.2.1]{Douc}. 

\begin{pro}
The class  of all irreducible Markov processes with initial distribution  is $\mathsf{InPL}_{\omega_1}$-definable by the formula
$\mathbf{Irr}:= \neg \mathbb{M}^{\langle{0}\rangle}_0 (p) \to\bigvee_{n\geq 1} \neg \mathbb{M}_0^{\langle{n}\rangle} (p)$.
\end{pro}
\begin{proof}
	Eazy to prove.
\end{proof}

\paragraph{Recurrence.}
Let $\langle{\Omega, \mathcal{A}, T}\rangle$ be a Markov process. The {\em expected number} of visits to a set $A\in \mathcal{A}$ starting from a state $w\in \Omega$, denoted by $U(w, A)$, is defined by
\[U(w, A):= \sum_{n=0}^{\infty} T^n(w, A).\] 
Note that for $w\in \Omega$, we have $U(w, A) \in \mathbb{N}\cup \{\infty \}$ is a measure on $\mathcal{A}$. This measure is called the {\em potential kernel} associated to $T$. A set $A\in \mathcal{A}$ is said to be {\em recurrent} if
$U(w, A)= \infty$  for all $w\in A$. Moreover,  $A\in \mathcal{A}$ is said to be {\em accessible} if for all $w\in \Omega$, there exists an integer $n\geq 1$ such that
$T^n(w, A)>0$. We say that a Markov process $\langle{\Omega, \mathcal{A}, T}\rangle$ is {\em recurrent} if it is irreducible and any accessible set is recurrent, e.g., see \cite[Definition 10.1.1]{Douc}.

\begin{pro}
	The class $\mathcal{R}$ of all recurrent Markov processes is $\mathsf{PL}_{\omega_1}$-definable by the conjunction of the formula $\mathbf{Irr}  = \neg \mathbb{M}^{\langle{0}\rangle}_0 (p) \to\bigvee_{n\geq 1} \neg \mathbb{M}_0^{\langle{n}\rangle} (p)$ and the formula
\[\mathbf{Rec}:= \big( \bigvee_{n\geq 1} \neg \mathbb{M}_0^{\langle{n}\rangle} (p)\big) \to \big( p\to \bigwedge_{r\in \mathbb{Q}^+} \; \bigvee_{k\in \mathbb{N}}\; \bigvee_{ r_0+\dots+r_{k}\geq r}\; \bigwedge_{m= 0}^{k} \mathbb{L}_{r_m}^{\langle{m}\rangle} (p)  \big). \]
\end{pro}
\begin{proof}
	Let $\mathfrak{P}= \langle{\Omega, \mathcal{A}, T}\rangle$ be an irreducible Markov process.  It is enough to show that
$\mathfrak{P}\vDash \mathbf{Rec}$ if and only if $\mathfrak{P}$ is recurrent.
First suppose that  $\mathfrak{P}$ is recurrent, i.e. for any $A\in \mathcal{A}$, if $A$ is accessible then $A$ is recurrent. Let $ \mathfrak{M}$ be a model based on $\mathfrak{P}$ and  $[\![p]\!]_{\mathfrak{M}} = A$ for some $A\in \mathcal{A}$. Then, we have the following facts:
 \begin{itemize}
 \item  $A$ is accessible if and only if $\mathfrak{M},w\vDash \bigvee_{n\geq 1} \neg \mathbb{M}_0^{\langle{n}\rangle} (p)$ for each  $w\in \Omega$.
 \item $A$ is recurrent  if and only if  for all $w\in A$ and $r\in \mathbb{Q}^+$, there exists $k\in \mathbb{N}$ such that $\sum_{m=0}^{k} T^m(w, A)\geq r$. Equivalently, for all $w\in \Omega$, if $w\in A$, then for all $r\in \mathbb{Q}^+$, there are $k\in \mathbb{N}$ and
  $r_0, \dots, r_{k} \in \mathbb{Q} \cap [0,1]$ such that $r_0+\dots+r_{k}\geq r$ and for all  $0\leq m\leq k$, we have  $T^m(w, A) \geq r_m$. This means that for all $w\in \Omega$,
  \[\mathfrak{M},w\vDash p\to   \bigwedge_{r\in \mathbb{Q}^+} \; \bigvee_{k\in \mathbb{N}}\; \bigvee_{r_0+\dots+r_{k}\geq r}\; \bigwedge_{m= 0}^{k} \mathbb{L}_{r_m}^{\langle{m}\rangle} (p).\]
 \end{itemize}
So it is easy to see that $\mathfrak{P}\vDash \mathbf{Rec}$.  The converse direction can be proved similarly.
\end{proof}

\appendix
\section{Appendix} \label{Appendix}
We here first review some basic notions from measurable spaces and measures that can be found in any textbook on measure theory and probability such as \cite{Bill, Dudley}. Then, we recall a few concepts from  Polish spaces. For more details about these concepts, see, e.g., \cite{Sri}.
\subsection{Measurable spaces and measures}
Let $\mathcal{A}$ be an (Boolean) algebra on a non-empty set $\Omega$. $\mathcal{A}$ is a {\em $\sigma$-algebra} if it is closed under countable unions. Then the pair $\langle{\Omega, \mathcal{A}}\rangle$ is called a {\em measurable space} and the elements of $\mathcal{A}$ {\em measurable sets} or {\em events}.
For a collection $A$ of subsets of $\Omega$, the {\em $\sigma$-algebra generated by} $A$, denoted $\sigma(A)$, is the smallest $\sigma$-algebra containing $A$.

Given a measurable space $\langle{\Omega, \mathcal{A}}\rangle$, a function $f : \Omega \to \Omega$ is {\em measurable} on $\mathcal{A}$, or simply {\em measurable}, if $f^{-1}(A) \in \mathcal{A}$ for all $A \in \mathcal{A}$. Now, let $\mathbb{R}^+=\{r\in \mathbb{R}\;|\; r \geq 0\}$. A set function $\mu: \mathcal{A}\to \mathbb{R}^+$ is {\em finitely additive} if $\mu (A_1 \cup A_2)= \mu (A_1)+ \mu(A_2)$ whenever $A_1$ and $A_2$ are disjoint members of $\mathcal{A}$. We say that $\mu$ is
{\em countably additive} (or {\em $\sigma$-additive}) if for any $\{A_n\}_{n\in \mathbb{N}}$ of pairwise disjoint members of $\mathcal{A}$, we have $\mu (\bigcup_{n\in \mathbb{N}}\ A_n)= \sum_{n \in \mathbb{N}} \mu (A_n)$.
A countably additive set function $\mu: \mathcal{A}\to \mathbb{R}^+$ is called a {\em measure} on $ \mathcal{A}$ if $\mu (\emptyset)=0$.

\begin{fact} [\cite{Bill}, Theorem 10.2] \label{contin}
	Let $\mathcal{A}$ be an algebra on $\Omega$ and $\mu: \mathcal{A}\to \mathbb{R}^+$ be finitely additive. Then the following are equivalent:
	\begin{itemize}
		\item[1.] $\mu$ is countably subadditive, i.e. $\mu (\bigcup_{n\in \mathbb{N}}\ A_n) \leq \sum_{n \in \mathbb{N}} \mu (A_n)$ for any countable family $A_n$ such that $\bigcup_{n \in \mathbb{N}} A_n\in \mathcal{A}$;
		\item[2.] $\mu$ is continuous from above at the empty set, i.e. $\mu(\bigcap_{n \in \mathbb{N}} A_n)= 0$ for any countable chain $A_0 \supseteq A_1 \supseteq \dots$ such that $\bigcap_{n \in \mathbb{N}} A_n = \emptyset$.
	\end{itemize}
\end{fact}

\begin{fact} [\cite{Bill}, Theorem 11.3] \label{extension}
	Let $\mathcal{A}$ be an algebra on $\Omega$ and $\mu: \mathcal{A}\to \mathbb{R}^+$ be finitely additive and countably subadditive. Then $\mu$ extends uniquely to a measure on $\sigma(\mathcal{A})$.
\end{fact}

We call that a measure $\mu: \mathcal{A}\to \mathbb{R}^+$ is a {\em probability measure} if in addition $\mu (\Omega)= 1$. In this case, the triple $\langle{\Omega, \mathcal{A}, \mu}\rangle$ is called a {\em probability space}.
\subsection{Polish spaces}

A topological space is said to be a {\em Polish space} if it is homeomorphic to a complete separable metric space. A subset of a topological space (metric space) is called a {\em $G_\delta$ set} if it can be represented as a countable intersection of open sets.

\begin{fact}[\cite{Sri}, Theorem 2.2.1] \label{G}
	Every $G_\delta$ subset of a Polish space is Polish.
\end{fact}

A {\em Hausdorff space} is a topological space with a separation property: any two distinct points can be separated by two disjoint open sets. A {\em locally compact space} is a topological space in which every point has a compact neighborhood.
\begin{fact}[\cite{Sri}, Corollary 2.3.32] \label{locally}
	Every second countable locally compact Hausdorff space is a Polish space.
\end{fact}
Let $\langle{X, \tau}\rangle$ be a topological space. The $\sigma$-algebra generated by the open sets of $\tau$, denoted $\mathcal{B}(\tau)$, is called the {\em Borel $\sigma$-algebra} and its measurable sets are called {\em Borel sets}. A measurable space $\langle{\Omega, \mathcal{A}}\rangle$ is called a {\em standard Borel space} if
$\mathcal{A}$ is the Borel $\sigma$-algebra generated by a Polish topology on $\Omega$.
\end{document}